\documentclass[a4paper,oneside]{amsart}

\RequirePackage{amsmath}
\RequirePackage{bm}
\RequirePackage{amssymb}
\RequirePackage{upref}
\RequirePackage{amsthm}
\RequirePackage{enumerate}
\RequirePackage{pb-diagram}
\RequirePackage{amsfonts}
\RequirePackage[mathscr]{eucal}
\RequirePackage{verbatim}
\RequirePackage{xr}
\RequirePackage{graphicx}
\RequirePackage[pdftex]{floatflt}
\usepackage{calc}
\usepackage{xspace}

\newcommand{\cf}{cf.\@\xspace}
\newcommand{\resp}{resp.\@\xspace}

\newcommand{\aev}{a.e.\@\xspace}

\newcommand{\wlogc}{w.l.o.g.\@\xspace}

\newcommand{\al}{\alpha}
\newcommand{\bet}{\beta}
\newcommand{\ga}{\gamma}
\newcommand{\de}{\delta }
\newcommand{\e}{\epsilon}

\newcommand{\f}{\varphi}
\newcommand{\h}{\eta}

\newcommand{\ka}{\kappa}
\newcommand{\lam}{\lambda}
\newcommand{\m}{\mu}
\newcommand{\n}{\nu}
\newcommand{\om}{\omega}

\newcommand{\s}{\sigma}
\newcommand{\x}{\xi}

\newcommand{\C}{\varGamma}
\newcommand{\D}{\varDelta}

\newcommand{\di}[1]{#1\nobreakdash-\hspace{0pt}dimensional}

\newcommand{\fu}[3]{#1\hspace{0pt}_{|_{#2_#3}}}
\newcommand{\fv}[2]{#1\hspace{0pt}_{|_{#2}}}

\newcommand{\so}{{\mc S_0}}

\newcommand{\const}{\tup{const}}
\newcommand{\ndash}{\nobreakdash--}

\newcommand{\msp[1]}[1]{\mspace{#1mu}}

\newcommand{\R}[1][n+1]{{\protect\mathbb R}^{#1}}

\newcommand{\N}{{\protect\mathbb N}}

\newcommand{\eR}{\stackrel{\lower1ex \hbox{\rule{6.5pt}{0.5pt}}}{\msp[3]\R[]}}
\newcommand{\eN}{\stackrel{\lower1ex \hbox{\rule{6.5pt}{0.5pt}}}{\msp[1]\N}}
\newcommand{\eO}{\stackrel{\lower1ex
\hbox{\rule{6pt}{0.5pt}}}{\msc O}}

\DeclareMathOperator{\graph}{graph}

\newcommand\ra{\rightarrow}

\newcommand{\ua}{\uparrow}
\newcommand{\da}{\downarrow}

\newcommand\pde[2]{\frac {\partial#1}{\partial#2}}
 
\newcommand{\un}{\infty}
\newcommand{\A}{\forall}

\newcommand{\set}[2]{\{\,#1\colon #2\,\}}
\newcommand{\uu}{\cup}

\newcommand{\uuu}{\bigcup}

\newcommand{\uud}{ \stackrel{\lower 1ex \hbox {.}}{\uu}}
\newcommand{\uuud}[1]{ \stackrel{\lower 1ex \hbox {.}}{\uuu_{#1}}}
\newcommand\su{\subset}

\newcommand{\sminus}[1][28]{\raise 0.#1ex\hbox{$\scriptstyle\setminus$}}

\newcommand{\wed}{\wedge}

\newcommand\ti{\times }

\newcommand{\abs}[1]{\lvert#1\rvert}

\newcommand{\norm}[1]{\lVert#1\rVert}

\newcommand{\nnorm}[1]{| \mspace{-2mu} |\mspace{-2mu}|#1| \mspace{-2mu}
|\mspace{-2mu}|}
\newcommand{\spd}[2]{\protect\langle #1,#2\protect\rangle}

\newcommand\ch[3]{\varGamma_{#1#2}^#3}
\newcommand\cha[3]{{\bar\varGamma}_{#1#2}^#3}

\newcommand{\riem}[4]{R_{#1#2#3#4}}
\newcommand{\riema}[4]{{\bar R}_{#1#2#3#4}}

\newcommand{\tit}{\textit}

\newcommand{\tup}{\textup}

\newcommand{\mc}{\protect\mathcal}
\newcommand{\msc}{\protect\mathscr}

\providecommand{\bysame}{\makebox[3em]{\hrulefill}\thinspace}

\newcommand{\ci}{\cite}

\newcommand{\cq}[1]{\glqq{#1}\grqq\,}

\newcommand{\bt}{\begin{thm}}
\newcommand{\bl}{\begin{lem}}
\newcommand{\bc}{\begin{cor}}
\newcommand{\bd}{\begin{definition}}
\newcommand{\bpp}{\begin{prop}}
\newcommand{\br}{\begin{rem}}
\newcommand{\bn}{\begin{note}}
\newcommand{\be}{\begin{ex}}
\newcommand{\bes}{\begin{exs}}
\newcommand{\bb}{\begin{example}}
\newcommand{\bbs}{\begin{examples}}
\newcommand{\ba}{\begin{axiom}}

\newcommand{\et}{\end{thm}}
\newcommand{\el}{\end{lem}}
\newcommand{\ec}{\end{cor}}
\newcommand{\ed}{\end{definition}}
\newcommand{\epp}{\end{prop}}
\newcommand{\er}{\end{rem}}
\newcommand{\en}{\end{note}}
\newcommand{\ee}{\end{ex}}
\newcommand{\ees}{\end{exs}}
\newcommand{\eb}{\end{example}}
\newcommand{\ebs}{\end{examples}}
\newcommand{\ea}{\end{axiom}}

\newcommand{\bp}{\begin{proof}}
\newcommand{\ep}{\end{proof}}
\newcommand{\eps}{\renewcommand{\qed}{}\end{proof}}

\newcommand{\bal}{\begin{align}}

\newcommand{\bi}[1][1.]{\begin{enumerate}[\upshape #1]}
\newcommand{\bia}[1][(1)]{\begin{enumerate}[\upshape #1]}
\newcommand{\bin}[1][1]{\begin{enumerate}[\upshape\bfseries #1]}
\newcommand{\bir}[1][(i)]{\begin{enumerate}[\upshape #1]}
\newcommand{\bic}[1][(i)]{\begin{enumerate}[\upshape\hspace{2\cma}#1]}
\newcommand{\bis}[2][1.]{\begin{enumerate}[\upshape\hspace{#2\parindent}#1]}
\newcommand{\ei}{\end{enumerate}}

\newcommand\ndots{\raise 0.47ex \hbox {,}\hskip0.06em\cdots %
     \raise 0.47ex \hbox {,}\hskip0.06em} 

\newcommand{\q}{\quad}
\newcommand{\qq}{\qquad}

\newcommand{\hp}{\hphantom}

\newcommand\nd{\noindent}

\newskip\Csmallskipamount                                                
\Csmallskipamount=\smallskipamount
\newskip\Cmedskipamount
\Cmedskipamount=\medskipamount
\newskip\Cbigskipamount
\Cbigskipamount=\bigskipamount

\newcommand\cvs{\vspace\Csmallskipamount}   
\newcommand\cvm{\vspace\Cmedskipamount}

\newskip\csa
\csa=\smallskipamount

\newskip\cma
\cma=\medskipamount

\newskip\cba
\cba=\bigskipamount

\newdimen\spt
\newcommand\citem{\cvs\advance\itemno by
1{(\romannumeral\the\itemno})\hskip3pt}
\newcommand{\bitem}{\cvm\nd\advance\itemno by
1{\bf\the\itemno}\hspace{\cma}}

\newcount\itemno
\newcommand{\las}[1]{\label{S:#1}}

\newcommand{\lae}[1]{\label{E:#1}}
\newcommand{\lat}[1]{\label{T:#1}}
\newcommand{\lal}[1]{\label{L:#1}}
\newcommand{\lad}[1]{\label{D:#1}}
\newcommand{\lac}[1]{\label{C:#1}}

\newcommand{\rs}[1]{Section~\ref{S:#1}}

\newcommand{\rt}[1]{Theorem~\ref{T:#1}}
\newcommand{\rl}[1]{Lemma~\ref{L:#1}}
\newcommand{\rd}[1]{Definition~\ref{D:#1}}
\newcommand{\rc}[1]{Corollary~\ref{C:#1}}

\newcommand{\re}[1]{\eqref{E:#1}}

\newskip\thmskip
\thmskip=\parindent

\newskip\hsk
\newenvironment{hinw}{\labelsep=0pt\begin{list}{}{\labelsep=0pt\itemindent=0pt\labelwidth=0pt\leftmargin=\parindent\rightmargin=0pt\partopsep=\cba}%
\item\it\nopagebreak\nopagebreak}%
{\end{list}}

\newcommand\bh{\begin{hinw}}
\newcommand{\eh}{\end{hinw}}

\newtheoremstyle{normal}
  {\cba}
  {\cba}
  {}
  {\thmskip}
  {\bfseries}
  {.}
  {\hsk}
  {}

\newtheoremstyle{abschnitt}
  {\cba}
  {\cba}
  {}
  {\thmskip}
  {\bfseries}
  {.}
  {\hsk}
  {}

\newtheoremstyle{italic}
  {\cba}
  {\cba}
  {\itshape}
  {\thmskip}
  {\bfseries}
  {.}
  {\hsk}
  {}

\newtheoremstyle{aufgaben}
  {\cba}
  {\cba}
  {}
  {}
  {\normalsize\bfseries}
  {.}
  {\hsk}
  {}

\newtheoremstyle{break}
  {\cba}
  {\cba}
  {\itshape}
  {}
  {\bfseries}
  {.}
  {\newline}
  {}

\swapnumbers
\theoremstyle{italic}
\newtheorem{thm}[subsection]{Theorem}
\newtheorem{lem}[subsection]{Lemma}
\newtheorem{prop}[subsection]{Proposition}
\newtheorem{cor}[subsection]{Corollary}

\theoremstyle{normal}
\newtheorem{rem}[subsection]{Remark}
\newtheorem{definition}[subsection]{Definition}
\newtheorem{example}[subsection]{Example}
\newtheorem{examples}[subsection]{Examples}
\newtheorem{ex}[subsection]{Exercise}
\newtheorem{note}[subsection]{}
\newtheorem{axiom}[subsection]{Axiom}

\theoremstyle{aufgaben}
\newtheorem{exs}[subsection]{Exercises}

\swapnumbers

\numberwithin{equation}{section}
\numberwithin{figure}{section}

\newenvironment{textequation}[1][0.8]
{\begin{equation}
\begin{aligned}
\begin{minipage}{#1\linewidth}}
{\end{minipage}
\end{aligned}
\end{equation}
\ignorespacesafterend}

\newcommand{\btext}{\begin{textequation}}
\newcommand{\etext}{\end{textequation}}

\usepackage[german,english]{babel}
\usepackage{graphicx}
\RequirePackage{amsmath}
\RequirePackage{bm}
\RequirePackage{amssymb}
\RequirePackage{upref}
\RequirePackage{amsthm}
\RequirePackage{enumerate}
\RequirePackage{pb-diagram}
\RequirePackage{amsfonts}
\RequirePackage[mathscr]{eucal}
\RequirePackage{verbatim}
\RequirePackage{xr}
\RequirePackage{graphicx}
\usepackage{calc}
\usepackage{xspace}

\RequirePackage{upref}
\RequirePackage{amsthm}
\RequirePackage{enumerate}
\usepackage[mathscr]{eucal}

\usepackage{xr-hyper}

\usepackage{calc}

\newlength{\oddsidemarginlength}
\newlength{\topmarginlength}

\newcounter{numberoflines}
\newcounter{tempcc}
\oddsidemargin=\oddsidemarginlength
\evensidemargin=\oddsidemargin
\topmargin=\topmarginlength
\usepackage{hyperref}

\begin{document}

\flushbottom


\title[The inverse mean curvature flow in ARW spaces]{The inverse mean curvature
flow in ARW spaces---transition from big crunch to big bang}

\author{Claus Gerhardt}
\address{Ruprecht-Karls-Universit\"at, Institut f\"ur Angewandte Mathematik,
Im Neuenheimer Feld 294, 69120 Heidelberg, Germany}
\email{gerhardt@math.uni-heidelberg.de}
\urladdr{http://www.math.uni-heidelberg.de/studinfo/gerhardt/}
\thanks{This work has been supported by the Deutsche Forschungsgemeinschaft.}

%
\date{\today}
%


\begin{abstract}
We consider spacetimes $N$ satisfying some structural
conditions, which are still fairly general, and prove convergence results for the leaves
of an inverse mean curvature flow. 

Moreover, we define a new spacetime $\hat N$ by switching the light
cone and using reflection to define a new time function, such that the two
spacetimes $N$ and $\hat N$ can be pasted together to yield a smooth manifold
having a metric singularity, which, when viewed from the region $N$ is a big crunch,
and when viewed from $\hat N$ is a big bang.

The inverse mean curvature flows in $N$ \resp $\hat N$ correspond to each other
via reflection. Furthermore, the properly rescaled flow in $N$ has a natural smooth
extension of class $C^3$ across the singularity into $\hat N$. With respect to this
natural, globally defined diffeomorphism we speak of a transition from big crunch to big
bang.
\end{abstract}
\maketitle

\tableofcontents

\setcounter{section}{-1}
\section{Introduction}

In \cite{cg:imcf} we considered the inverse mean curvature flow (IMCF) in
cosmological spacetimes having a future mean curvature barrier and showed that the
IMCF  exists for all time and runs directly into the future singularity, if and
only if $N$ satisfies a \tit{strong volume decay condition}.

Apart from the fact that the leaves run straight into the future singularity no further
convergence results could be derived due to the weak assumptions on the spacetime.

\cvm
In the present paper we consider spacetimes $N$ satisfying some structural
conditions, which are still fairly general, and prove convergence results for the leaves
of the IMCF. 

Moreover, we define a new spacetime $\hat N$ by switching the light
cone and using reflection to define a new time function, such that the two
spacetimes $N$ and $\hat N$ can be pasted together to yield a smooth manifold
having a metric singularity, which, when viewed from the region $N$ is a big crunch,
and when viewed from $\hat N$ is a big bang.

The inverse mean curvature flows in $N$ \resp $\hat N$ correspond to each other
via reflection. Furthermore, the properly rescaled flow in $N$ has a natural smooth
extension of class $C^3$ across the singularity into $\hat N$. With respect to this
natural diffeomorphism we speak of a transition from big crunch to big bang.

\bd\lad{0.1}
A cosmological spacetime $N$, $\dim N=n+1$, is said to be \tit{asymptotically
Robertson-Walker} (ARW) with respect to the future, if a future end of $N$, $N_+$,
can be written as a product $N_+=[a,b)\times \so$, where $\so$ is a compact
Riemannian space, and there exists a future directed time function $\tau=x^0$ such
that the metric in $N_+$ can be written as
\begin{equation}\lae{1.19}
d\breve s^2=e^{2\tilde\psi}\{-{(dx^0})^2+\s_{ij}(x^0,x)dx^idx^j\},
\end{equation}
where  $\so$ corresponds to $x^0=a$, $\tilde\psi$ is of the form
\begin{equation}
\tilde\psi(x^0,x)=f(x^0)+\psi(x^0,x),
\end{equation}
and we assume that there exists a positive constant $c_0$ and a smooth
Riemannian metric $\bar\s_{ij}$ on $\so$ such that
\begin{equation}
\lim_{\tau\ra b}e^\psi=c_0\q\wed\q \lim_{\tau\ra b}\s_{ij}(\tau,x)=\bar\s_{ij}(x),
\end{equation}
and
\begin{equation}
\lim_{\tau\ra b}f(\tau)=-\un.
\end{equation}

\cvm
Without loss of generality we shall assume $c_0=1$. Then $N$ is ARW with
respect to the future, if the metric is close to the Robertson-Walker metric
\begin{equation}\lae{0.5}
d\bar s^2=e^{2f}\{-{dx^0}^2+\bar\s_{ij}(x)dx^idx^j\}
\end{equation}
near the singularity $\tau =b$. By \tit{close} we mean that the derivatives of arbitrary order with respect to space and time of the
conformal metric $e^{-2f}\breve g_{\al\bet}$ in \re{1.19} should converge  to the
corresponding derivatives of the conformal limit metric in \re{0.5} when $x^0$ tends
to $b$. We emphasize that in our terminology Robertson-Walker metric does not
imply that
$(\bar\s_{ij})$ is a metric of constant curvature, it is only the spatial metric of a
warped product.

\cvm
We assume, furthermore, that $f$ satisfies the following five conditions
\begin{equation}
-f'>0,
\end{equation}
there exists $\om\in\R[]$ such that
\begin{equation}\lae{0.7}
n+\om-2>0\q\wed\q \lim_{\tau\ra b}\abs{f'}^2e^{(n+\om-2)f}=m>0.
\end{equation}
Set $\tilde\ga =\frac12(n+\om-2)$, then there exists the limit
\begin{equation}\lae{0.8}
\lim_{\tau\ra b}(f''+\tilde\ga \abs{f'}^2)
\end{equation}
and
\begin{equation}\lae{0.9}
\abs{D^m_\tau(f''+\tilde\ga \abs{f'}^2)}\le c_m \abs{f'}^m\qq
\A\, m\ge 1,
\end{equation}
as well as
\begin{equation}\lae{0.10}
\abs{D_\tau^mf}\le c_m \abs{f'}^m\qq\A\, m\ge 1.
\end{equation}

\cvm
We call $N$ a \tit{normalized} ARW spacetime, if
\begin{equation}
\int_{\so}\sqrt{\det{\bar\s_{ij}}}=\abs{S^n}.
\end{equation}
\ed

\br
(i) If these assumptions are satisfied, then we shall show that the range of $\tau$ is
finite, hence, we may---and shall---assume \wlogc that $b=0$, i.e.,
\begin{equation}
a<\tau<0.
\end{equation}

\cvm
(ii) Any ARW spacetime can be normalized as one easily checks. For normalized ARW
spaces the constant
$m$ in \re{0.7} is defined uniquely and can be identified
with the mass of $N$, \cf \cite{cg:mass}.

\cvm
(iii) In view of the assumptions on $f$ the mean curvature of the coordinate slices
$M_\tau=\{x^0=\tau\}$ tends to $\un$, if $\tau$ goes to zero.

\cvm
(iv) ARW spaces satisfy a strong volume decay condition, \cf
\cite[Definition~0.1]{cg:imcf}.

\cvm
(v) Similarly one can define $N$ to be ARW with respect to the past. In this case the
singularity would lie in the past, correspond to $\tau=0$, and the mean curvature
of the coordinate slices would tend to $-\un$.
\er

We assume that $N$ satisfies the timelike convergence condition. Consider the future
end $N_+$ of $N$ and let $M_0\su N_+$ be a spacelike hypersurface with positive
mean curvature $\fv {\breve H}{M_0}>0$ with respect to the past directed normal
vector $\breve\nu$---we shall explain in \rs{2} why we use the symbols $\breve H$
and $\breve\nu$ and not the usual ones $H$ and $\nu$. Then, as we have proved in
\cite{cg:imcf}, the inverse mean curvature flow 
\begin{equation}
\dot x=-\breve H^{-1}\breve\nu
\end{equation}
with initial hypersurface $M_0$ exists for all time, is smooth, and runs straight
into the future singularity.

\cvm
If we express the flow hypersurfaces $M(t)$  as graphs over $\so$
\begin{equation}
M(t)=\graph u(t,\cdot),
\end{equation}
then our main results can be formulated as

\bt\lat{0.3}
\tup{(i)} Let $N$ satisfy the above assumptions, then the range of the time function
$x^0$ is finite, i.e., we may assume that $b=0$. Set
\begin{equation}
\tilde u=ue^{\ga t},
\end{equation}
where $\ga=\tfrac1n\tilde\ga$, then there are positive constants $c_1, c_2$ such
that
\begin{equation}
-c_2\le \tilde u\le -c_1<0,
\end{equation}
and $\tilde u$ converges in $C^\un(\so)$ to a smooth function, if $t$ goes to
infinity. We shall also denote the limit function  by $\tilde u$.

\cvm
\tup{(ii)} Let $\breve g_{ij}$ be the induced metric of the leaves $M(t)$, then the
rescaled metric
\begin{equation}
e^{\frac2n t} \breve g_{ij}
\end{equation}
converges in $C^\un(\so)$ to
\begin{equation}
(\tilde\ga m)^\frac1{\tilde\ga}(-\tilde u)^\frac2{\tilde\ga}\bar\s_{ij}.
\end{equation}

\cvm
\tup{(iii)} The leaves $M(t)$ get more umbilical, if $t$ tends to infinity, namely, there
holds
\begin{equation}
\breve H^{-1}\abs{\breve h^j_i-\tfrac1n \breve H\de^j_i}\le c\msp e^{-2\ga t}.
\end{equation}
In case $n+\om-4>0$, we even get a better estimate
\begin{equation}
\abs{\breve h^j_i-\tfrac1n \breve H\de^j_i}\le c\msp e^{-\frac1{2n}(n+\om-4) t}.
\end{equation}
\et

\cvm
For a description of the results related to the transition from big crunch to big bang
we refer to \rs{8}.

\section{Notations and definitions}\las{1}

The main objective of this section is to state the equations of Gau{\ss}, Codazzi,
and Weingarten for space-like hypersurfaces $M$ in a \di {(n+1)} Lorentzian
manifold
$N$.  Geometric quantities in $N$ will be denoted by
$(\bar g_{ \al \bet}),(\riema  \al \bet \ga \de)$, etc., and those in $M$ by $(g_{ij}), 
(\riem ijkl)$, etc.. Greek indices range from $0$ to $n$ and Latin from $1$ to $n$;
the summation convention is always used. Generic coordinate systems in $N$ resp.
$M$ will be denoted by $(x^ \al)$ resp. $(\x^i)$. Covariant differentiation will
simply be indicated by indices, only in case of possible ambiguity they will be
preceded by a semicolon, i.e., for a function $u$ in $N$, $(u_ \al)$ will be the
gradient and
$(u_{ \al \bet})$ the Hessian, but e.g., the covariant derivative of the curvature
tensor will be abbreviated by $\riema  \al \bet \ga{ \de;\e}$. We also point out that
\begin{equation}
\riema  \al \bet \ga{ \de;i}=\riema  \al \bet \ga{ \de;\e}x_i^\e
\end{equation}
with obvious generalizations to other quantities.

Let $M$ be a \tit{spacelike} hypersurface, i.e., the induced metric is Riemannian,
with a differentiable normal $\n$ which is time-like.

In local coordinates, $(x^ \al)$ and $(\x^i)$, the geometric quantities of the
space-like hypersurface $M$ are connected through the following equations
\begin{equation}\lae{1.2}
x_{ij}^ \al= h_{ij}\n^ \al
\end{equation}
the so-called \tit{Gau{\ss} formula}. Here, and also in the sequel, a covariant
derivative is always a \tit{full} tensor, i.e.

\begin{equation}
x_{ij}^ \al=x_{,ij}^ \al-\ch ijk x_k^ \al+ \cha  \bet \ga \al x_i^ \bet x_j^ \ga.
\end{equation}
The comma indicates ordinary partial derivatives.

In this implicit definition the \tit{second fundamental form} $(h_{ij})$ is taken
with respect to $\n$.

The second equation is the \tit{Weingarten equation}
\begin{equation}
\n_i^ \al=h_i^k x_k^ \al,
\end{equation}
where we remember that $\n_i^ \al$ is a full tensor.

Finally, we have the \tit{Codazzi equation}
\begin{equation}
h_{ij;k}-h_{ik;j}=\riema \al \bet \ga \de\n^ \al x_i^ \bet x_j^ \ga x_k^ \de
\end{equation}
and the \tit{Gau{\ss} equation}
\begin{equation}
\riem ijkl=- \{h_{ik}h_{jl}-h_{il}h_{jk}\} + \riema  \al \bet\ga \de x_i^ \al x_j^ \bet
x_k^ \ga x_l^ \de.
\end{equation}

Now, let us assume that $N$ is a globally hyperbolic Lorentzian manifold with a
\tit{compact} Cauchy surface. 
$N$ is then a topological product $I\times \mc S_0$, where $I$ is an open interval,
$\mc S_0$ is a compact Riemannian manifold, and there exists a Gaussian coordinate
system
$(x^ \al)$, such that the metric in $N$ has the form 
\begin{equation}\lae{1.7}
d\bar s_N^2=e^{2\psi}\{-{dx^0}^2+\s_{ij}(x^0,x)dx^idx^j\},
\end{equation}
where $\s_{ij}$ is a Riemannian metric, $\psi$ a function on $N$, and $x$ an
abbreviation for the spacelike components $(x^i)$. 
We also assume that
the coordinate system is \tit{future oriented}, i.e., the time coordinate $x^0$
increases on future directed curves. Hence, the \tit{contravariant} time-like
vector $(\x^ \al)=(1,0,\dotsc,0)$ is future directed as is its \tit{covariant} version
$(\x_ \al)=e^{2\psi}(-1,0,\dotsc,0)$.

Let $M=\graph \fv u\so$ be a space-like hypersurface
\begin{equation}
M=\set{(x^0,x)}{x^0=u(x),\,x\in\mc S_0},
\end{equation}
then the induced metric has the form
\begin{equation}
g_{ij}=e^{2\psi}\{-u_iu_j+\s_{ij}\}
\end{equation}
where $\s_{ij}$ is evaluated at $(u,x)$, and its inverse $(g^{ij})=(g_{ij})^{-1}$ can
be expressed as
\begin{equation}\lae{1.10}
g^{ij}=e^{-2\psi}\{\s^{ij}+\frac{u^i}{v}\frac{u^j}{v}\},
\end{equation}
where $(\s^{ij})=(\s_{ij})^{-1}$ and
\begin{equation}\lae{1.11}
\begin{aligned}
u^i&=\s^{ij}u_j\\
v^2&=1-\s^{ij}u_iu_j\equiv 1-\abs{Du}^2.
\end{aligned}
\end{equation}
Hence, $\graph u$ is space-like if and only if $\abs{Du}<1$.

The covariant form of a normal vector of a graph looks like
\begin{equation}
(\n_ \al)=\pm v^{-1}e^{\psi}(1, -u_i).
\end{equation}
and the contravariant version is
\begin{equation}
(\n^ \al)=\mp v^{-1}e^{-\psi}(1, u^i).
\end{equation}
Thus, we have
\br Let $M$ be space-like graph in a future oriented coordinate system. Then the
contravariant future directed normal vector has the form
\begin{equation}
(\n^ \al)=v^{-1}e^{-\psi}(1, u^i)
\end{equation}
and the past directed
\begin{equation}\lae{1.15}
(\n^ \al)=-v^{-1}e^{-\psi}(1, u^i).
\end{equation}
\er

In the Gau{\ss} formula \re{1.2} we are free to choose the future or past directed
normal, but we stipulate that we always use the past directed normal for reasons
that we have explained in \ci[Section 2]{cg:indiana}.

Look at the component $ \al=0$ in \re{1.2} and obtain in view of \re{1.15}

\begin{equation}\lae{1.16}
e^{-\psi}v^{-1}h_{ij}=-u_{ij}- \cha 000\mspace{1mu}u_iu_j- \cha 0j0
\mspace{1mu}u_i- \cha 0i0\mspace{1mu}u_j- \cha ij0.
\end{equation}
Here, the covariant derivatives are taken with respect to the induced metric of
$M$, and
\begin{equation}
-\cha ij0=e^{-\psi}\bar h_{ij},
\end{equation}
where $(\bar h_{ij})$ is the second fundamental form of the hypersurfaces
$\{x^0=\const\}$.

An easy calculation shows
\begin{equation}
\bar h_{ij}e^{-\psi}=-\tfrac{1}{2}\dot\s_{ij} -\dot\psi\s_{ij},
\end{equation}
where the dot indicates differentiation with respect to $x^0$.

\section{The evolution problem}\las{2}

When proving the convergence results for the inverse mean curvature flow, we shall
consider the flow hypersurfaces to be embedded in $N$ equipped with the conformal
metric
\begin{equation}\lae{2.1}
d\bar s^2=-(dx^0)^2+ \s_{ij}(x^0,x)dx^idx^j.
\end{equation}

Though, formally, we have a different ambient space we still denote it by the same
symbol $N$ and distinguish only the metrics $\breve g_{\al\bet}$ and $\bar
g_{\al\bet}$
\begin{equation}
\breve g_{\al\bet}=e^{2\tilde\psi}\bar g_{\al\bet}
\end{equation}
and the corresponding geometric quantities of the hypersurfaces $\breve h_{ij},
\breve g_{ij}, \breve\nu$ \resp $h_{ij}, g_{ij}, \nu$, etc., i.e., the notations of the
preceding section now apply  to the case when $N$ is equipped with the metric in
\re{2.1}.

\cvm
The second fundamental forms $\breve h^j_i$ and $h^j_i$ are related by
\begin{equation}
e^{\tilde\psi}\breve h^j_i=h^j_i+\tilde\psi_\al\nu^\al\de^j_i
\end{equation}
and, if we define $F$ by
\begin{equation}
F=e^{\tilde\psi} \breve H,
\end{equation}
then
\begin{equation}
F=H-n\tilde vf'+n\psi_\al\nu^\al,
\end{equation}
where
\begin{equation}
\tilde v=v^{-1},
\end{equation}
and the evolution equation can be written as
\begin{equation}\lae{2.7}
\dot x=-F^{-1}\nu,
\end{equation}
since
\begin{equation}
\breve\nu=e^{-\tilde\psi}\nu.
\end{equation}

\cvm
The flow exists for all time and is smooth.

\cvm
Next, we want to show how the metric, the second fundamental form, and the
normal vector of the hypersurfaces $M(t)$ evolve. All time derivatives are
\tit{total} derivatives. We refer to \ci{cg:indiana} for more general results and to
\ci[Section 3]{cg96}, where proofs are given in a Riemannian setting, but these
proofs are also valid in a Lorentzian environment.

\bl\lal{2.1}
The metric, the normal vector, and the second fundamental form of $M(t)$
satisfy the evolution equations
\begin{equation}\lae{2.9}
\dot g_{ij}=-2 F^{-1}h_{ij},
\end{equation}
\begin{equation}\lae{2.10}
\dot \n=\nabla_M(-F^{-1})=g^{ij}(-F^{-1})_i x_j,
\end{equation}
and
\begin{equation}\lae{2.11}
\dot h_i^j=(-F^{-1})_i^j+F^{-1} h_i^k h_k^j + F^{-1} \riema
\al\bet\ga\de\n^\al x_i^\bet \n^\ga x_k^\de g^{kj}
\end{equation}
\begin{equation}
\dot h_{ij}=(-F^{-1})_{ij}-F^{-1} h_i^k h_{kj}+ F^{-1} \riema
\al\bet\ga\de\n^\al x_i^\bet \n^\ga x_j^\de.
\end{equation}
\el

Since the initial hypersurface is a graph over $\so$, we can write
\begin{equation}
M(t)=\graph\fu{u(t)}S0\q \A\,t\in I,
\end{equation}
where $u$ is defined in the cylinder $\R[]_+\times \so$. We then deduce from
\re{2.7}, looking at the component $\al=0$, that $u$ satisfies a parabolic
equation of the form
\begin{equation}\lae{2.14}
\dot u=\frac{\tilde v}F,
\end{equation}
where we  use the notations in \rs{1}, and where we emphasize that the time
derivative is a total derivative, i.e.

\begin{equation}
\dot u=\pde ut+u_i\dot x^i.
\end{equation}

Since the past directed normal can be expressed as
\begin{equation}
(\n^\al)=-e^{-\psi}v^{-1}(1,u^i),
\end{equation}
we conclude from \re{2.14}
\begin{equation}\lae{2.17}
\pde ut=\frac vF.
\end{equation}

Sometimes, we need a Riemannian reference metric, e.g., if we want to estimate
tensors. Since the Lorentzian metric can be expressed as
\begin{equation}
\bar g_{\al\bet}dx^\al dx^\bet=-(dx^0)^2+\s_{ij}dx^i dx^j,
\end{equation}
we define a Riemannian reference metric $(\tilde g_{\al\bet})$ by
\begin{equation}
\tilde g_{\al\bet}dx^\al dx^\bet=(dx^0)^2+\s_{ij}dx^i dx^j
\end{equation}
and we abbreviate the corresponding norm of a vectorfield $\h$ by
\begin{equation}
\nnorm \h=(\tilde g_{\al\bet}\h^\al\h^\bet)^{1/2},
\end{equation}
with similar notations for higher order tensors.

\section{Lower order estimates}\las{3}

We first draw a few immediate conclusions from our assumptions on $f$.

\bl\lal{3.1}
Let $f\in C^2([a,b))$ satisfy the conditions
\begin{equation}
\lim_{\tau\ra b}f(\tau)=-\un
\end{equation}
and
\begin{equation}\lae{3.2}
\lim_{\tau\ra b}\abs{f'}^2e^{2\tilde \ga f}=m,
\end{equation}
where $\tilde \ga, m$ are positive, then $b$ is finite.
\el

\bp
From \re{3.2} we deduce that $f'$ tends to $-\un$ and
\begin{equation}
\lim(-f'e^{\tilde \ga f})=\sqrt m.
\end{equation}

Moreover,
\begin{equation}
e^{\tilde\ga\tau}-e^{\tilde\ga \tau_0}=\int_{\tau_0}^\tau \tilde\ga f' e^{\tilde\ga
f}\le -\tilde\ga \tfrac {\sqrt{m}}2(\tau-\tau_0),
\end{equation}
if $\tau_0$ is close to $b$ in the topology of $\eR$ and $\tau>\tau_0$. Hence $b$
has to be finite.
\ep

\bc
We may---and shall---therefore assume that $b=0$, i.e., the time intervall $I$ is given
by $I=[a,0)$.
\ec

A simple application of de L'Hospital's rule then yields
\begin{equation}\lae{3.5}
\lim_{\tau\ra 0}\frac{e^{\tilde \ga f}}\tau=-\tilde\ga \sqrt m
\end{equation}

\cvm
From this relation and \re{0.8} we conclude

\bl
There holds
\begin{equation}\lae{3.6}
f' e^{\tilde\ga f}+\sqrt m\sim c\tau^2,
\end{equation}
where $c$ is a constant, and where the relation
\begin{equation}
\f\sim c\tau^2
\end{equation}
means 
\begin{equation}
\lim_{\tau\ra 0}\frac{\f(\tau)}{\tau^2}=c.
\end{equation}
\el

\bp
Applying de L'Hospital's rule we get
\begin{equation}
\lim \frac{f' e^{\tilde\ga f}+\sqrt m}{\tfrac12 \tau^2}=\lim \frac{(f''+\tilde\ga
\abs{f'}^2) e^{\tilde\ga f}}\tau =-c \tilde\ga \sqrt m.
\end{equation}
\ep

\bl\lal{3.4}
The asymptotic relation
\begin{equation}
\tilde\ga f'\tau-1\sim c\tau^2
\end{equation}
is valid.
\el

\bp
The relation \re{3.6} yields
\begin{equation}
\tilde\ga f'\tau e^{\tilde\ga f}+\sqrt m\tilde\ga \tau \sim c_1\tau^3,
\end{equation}
or equivalently,
\begin{equation}
(\tilde\ga f'\tau -1)e^{\tilde\ga f}+\sqrt m\tilde\ga \tau +e^{\tilde\ga f}\sim
c_1\tau^3.
\end{equation}

Dividing by $\tau^3$ and applying de L'Hospital's rule we infer
\begin{equation}
\lim \frac{\tilde\ga f'\tau -1}{\tau^2}\cdot \lim \frac{e^{\tilde\ga f}}\tau +\lim
\frac{\sqrt m\tilde \ga +\tilde\ga f' e^{\tilde\ga f}}{3\tau^2}=c_1,
\end{equation}
hence the result in view of \re{3.5} and \re{3.6}.
\ep

After these preliminary results we now want to prove that there are positive
constants $c_1,c_2$ such that
\begin{equation}
-c_1\le \tilde u\equiv u e^{\ga t}\le -c_2 <0\qq\A\,t\in\R[]_+,
\end{equation}
where $u$ is the solution of the scalar version of the inverse mean curvature flow,
i.e., $u$ is the solution of equation \re{2.14}.

\cvm
We shall proceed in two steps, first we shall derive
\begin{equation}\lae{3.15}
\abs{u e^{\lam t}}\le c(\lam)\qq\A\,0<\lam<\ga,
\end{equation}
and then the final result in the limiting case $\lam=\ga$.

\cvm
This procedure will also be typical for higher order estimates in the next sections.

\bl
For any $0<\lam<\ga$, there exists a constant $c(\lam)$ such that the estimate
\re{3.15} is valid.
\el

\bp
Define $\f=\f(t)$ by
\begin{equation}
\f(t)=\inf_{x\in\so}u(t,x).
\end{equation}
Then $\f$ is Lipschitz continuous and
\begin{equation}
\dot\f(t)=\pde ut(t,x_t)\qq \tup{for \aev}\,t,
\end{equation}
where $x_t\in\so$ is such that the infimum of $u(t,\cdot)$ is attained. This is a
well-known result, for a simple proof see e.g., \cite[Lemma 3.2]{cg:imcf}.

Let
\begin{equation}
w=\log(-\f) +\lam t,
\end{equation}
then, for \aev $t$, we have
\begin{equation}
\dot w=\f^{-1}\dot\f+\lam=u^{-1}\pde ut+\lam,
\end{equation} 
where $u$ is evaluated at $(t,x_t)$. In $x_t$ $u(t,\cdot)$ attains its infimum, i.e.,
$Du=0$ and $-\D u\le 0$.

From the parabolic equation \re{2.17}, we obtain in $x_t$
\begin{equation}
\pde ut=\frac1F=\frac1{H-nf'-n\dot\psi}.
\end{equation}
The mean curvature $H$ can be expressed as
\begin{equation}
H=-\D u+\bar H=-\D u +\s^{ij}\bar h_{ij}=-\D u-\tfrac12 \s^{ij}\dot\s_{ij}.
\end{equation}
Thus we deduce
\begin{equation}
\pde ut\ge \frac1{-nf'-n\dot\psi-\tfrac12\s^{ij}\dot\s_{ij}}
\end{equation}
and
\begin{equation}\lae{3.23}
\begin{aligned}
\dot w&\le
\frac1{-nf'u-(n\dot\psi-\tfrac12\s^{ij}\dot\s_{ij})u}+\lam\\[\cma]
&=\frac{1-nf'u\lam-(n\dot\psi-\tfrac12\s^{ij}\dot\s_{ij})\lam
u}{-nf'u-(n\dot\psi-\tfrac12\s^{ij}
\dot\s_{ij})u}.
\end{aligned}
\end{equation}

\cvm
Now, we observe that the argument of $f'$ is $u$ and
\begin{equation}
\lim_{t\ra\un}\inf_{x\in\so}u(t,x)=0,
\end{equation}
\cf \cite[Lemma 3.1]{cg:imcf}. Hence
\begin{equation}
\lim_{t\ra \un}f'u=\tilde\ga^{-1},
\end{equation}
in view of \rl{3.4}, and we infer that the right-hand side of inequality \re{3.23} is
negative for large $t$, $t\ge t_\lam$, and therefore
\begin{equation}
w\le w(t_\lam)\qq\A\,t\ge t_\lam,
\end{equation}
or equivalently,
\begin{equation}
-u e^{\lam t}\le c(\lam)\qq\A\,t\in\R[]_+.\qedhere
\end{equation}
\ep

\bt\lat{3.6}
Let $u$ be a solution of the evolution equation \re{2.14}, where $f$ satisfies the
assumptions \re{0.7} and \re{0.8}, then there are positive constants $c_1, c_2$ such
that
\begin{equation}\lae{3.28}
-c_1\le \tilde u\equiv ue^{\ga t}\le -c_2<0.
\end{equation}
\et

\bp
We only prove the estimate from above. Define
\begin{equation}
\f(t)=\sup_{x\in\so}u(t,x)
\end{equation}
and
\begin{equation}
w=\log(-\f)+\ga t.
\end{equation}

Arguing similar as in the proof of the previous lemma, we obtain for \aev $t$
\begin{equation}
\dot w\ge \frac{1-nf'u\ga-(n\dot\psi-\tfrac12\s^{ij}\dot\s_{ij})\ga
u}{-nf'u-(n\dot\psi-\tfrac12\s^{ij}
\dot\s_{ij})u}.
\end{equation}

Since $\tilde\ga= n\ga$, we deduce from \rl{3.4} that the right-hand side can be
estimated from below by $c \msp u$, i.e.,
\begin{equation}
\dot w\ge c\msp u\ge -c\msp c_\lam e^{-\lam t}
\end{equation}
for any $0<\lam <\ga$. Hence $w$ is bounded from below, or equivalently,
\begin{equation}
\tilde u\le -c_2<0.\qedhere
\end{equation}
\ep

\bc\lac{3.7}
For any $k\in\N^*$ there exists $c_k$ such that
\begin{equation}
\abs{f^{(k)}}\le c_k e^{k\ga t},
\end{equation}
where $f^{(k)}$ is evaluated at $u$.
\ec

\bp
In view of the assumption \re{0.10} there holds
\begin{equation}
\abs{f^{(k)}}\le c_k\abs{f'}^k=c_k\abs{f'}^k u^k \tilde u^{-k} e^{k\ga t}.
\end{equation}
Then use \rl{3.4} and the preceding theorem.
\ep

\section{$C^1$-estimates}\las{4}

We want to prove estimates for $\tilde v$ and $\norm{D\tilde u}$, where we recall
that
\begin{equation}
\tilde u= u e^{\ga t}.
\end{equation}

Our final goal is to show that $\norm{D\tilde u}$ is uniformly bounded, but this
estimate has to be deferred to \rs{5}. At the moment we only prove an exponential
decay for any $0<\lam <\ga$, i.e., we shall estimate $\norm{Du}\msp e^{\lam t}$.

\cvm
The starting point is the evolution equation satisfied by $\tilde v$.

\bl[Evolution of $\tilde v$]
Consider the flow \re{2.7}. Then $\tilde v$ satisfies the evolution equation
\begin{equation}\lae{4.2}
\begin{aligned}
\dot{\tilde v}&-F^{-2}\D \tilde v=
-F^{-2}\norm{A}^2\tilde v+F^{-2}\bar R_{\al\bet}\nu^\al
x^\bet_iu^i\\[\cma] &-F^{-2}(2H-nf'\tilde v+n\psi_\al \nu^\al)
\h_{\al\bet}\nu^\al\nu^\bet\\[\cma] &-F^{-2}(\h_{\al\bet\ga}\nu^\al x^\bet_i
x^\ga_j g^{ij}+\h_{\al\bet} x^\al_i x^\bet_j h^{ij})\\[\cma]
&-F^{-2}(-nf''\norm{Du}^2\tilde v-nf'\tilde v_ku^k
+n\psi_{\al\bet}\nu^\al x^\bet_iu^i+n \psi_\al x^\al_k h^k_iu^i),
\end{aligned}
\end{equation}
where $\h=(\h_\al)=(-1,0,\ldots,0)$ is a covariant unit vectorfield.
\el

\bp
We have
\begin{equation}
\tilde v=\h_\al\nu^\al.
\end{equation}
Let $(\x^i)$ be local coordinates for $M(t)$;  differentiating $\tilde v$ covariantly
we deduce
\begin{equation}\lae{4.4}
\tilde v_i=\h_{\al\bet} x^\bet_i\nu^\al+\h_\al \nu^\al_i,
\end{equation}
and
\begin{equation}
\tilde v_{ij}=\h_{\al\bet\ga} x^\bet_i x^\ga_j \nu^\al +\h_{\al\bet}\nu^\al_j
x^\bet_i+\h_{\al\bet}\nu^\al\nu^\bet h_{ij}+\h_\al \nu^\al_{ij}.
\end{equation}

The time derivative of $\tilde v$ is equal to
\begin{equation}
\begin{aligned}
\dot{\tilde v}&=\h_{\al\bet}\nu^\al\dot x^\bet +\h_\al \dot\nu^\al\\[\cma]
&=-\h_{\al\bet}\nu^\al\nu^\bet F^{-1} +F^{-2}\h_\al F^k x^\al_k.
\end{aligned}
\end{equation}

\cvm
From these relations the evolution equation for $\tilde v$ follows immediately with
the help of the Weingarten and Codazzi equations, the Gau{\ss} formula, and the
definition of $F$.
\ep

\bl\lal{4.2}
The following estimates are valid
\begin{equation}
\abs{\h_{\al\bet}\nu^\al\nu^\bet}\le c \tilde v^2 \nnorm{\h_{\al\bet}},
\end{equation}
\begin{equation}
\abs{\h_{\al\bet\ga}\nu^\al x^\bet_i x^\ga_j g^{ij}}\le c \tilde v^3
\nnorm{\h_{\al\bet\ga}},
\end{equation}
\begin{equation}
\abs{\h_{\al\bet}\nu^\al x^\bet_ku^k}\le c \nnorm{\h_{\al\bet}}\tilde v^3,
\end{equation}
\begin{equation}
\abs{\psi_\al x^\al_k h^k_iu^i}\le c \nnorm{D\psi}\msp \norm A\tilde v^2,
\end{equation}
\begin{equation}
\abs{\h_{\al\bet} x^\al_i x^\bet_j h^{ij}}\le c \nnorm{\h_{\al\bet}}\msp \norm
A\tilde v^2,
\end{equation}
and
\begin{equation}
\begin{aligned}
\abs{\bar R_{\al\bet}\nu^\al x^\bet_k u^k}\le c \tilde v^3 \abs{\bar R_{0k}\check
u^k}&+c \tilde v \abs{\bar R_{00}}\norm{Du}^2\\[\cma]
&+c\tilde v^3 \abs{\bar R_{ij}\check u^i \check u^j},
\end{aligned}
\end{equation}
where
\begin{equation}
\check u^i=\s^{ij}u_j.
\end{equation}
\el

\bp
Easy exercise.
\eps

We can now prove that $\tilde v$ is uniformly bounded.

\bl\lal{4.3}
The quantity $\tilde v$ is uniformly bounded
\begin{equation}
\tilde v\le c.
\end{equation}
\el

\bp
For large $T$, $0<T<\un$, assume that
\begin{equation}
\sup_{[0,T]}\sup_{M(t)}\tilde v=\tilde v(t_0,x_0).
\end{equation}
Applying the maximum principle we shall deduce that either $\tilde v \le 2$ or that
$t_0$ is a priori bounded 
\begin{equation}
t_0\le T_0.
\end{equation}

\cvm
In $(t_0,x_0)$ the left-hand side of equation \re{4.2} is non-negative, assuming
$t_0\ne 0$. Multiplying the resulting inequality by $F^2$ and using the estimates in
\rl{4.2} we conclude
\begin{equation}
\begin{aligned}
0\le -\norm A^2 \tilde v-n f''\norm{Du}^2\tilde v+c (1+\abs{f'})\tilde v^3 +c \norm
A \tilde v^2.
\end{aligned}
\end{equation}

If $\tilde v\ge 2$, then
\begin{equation}
\norm{Du}^2\ge \e_0 \tilde v^2
\end{equation}
with a positive constant $\e_0$, and if $t_0$ would be large, then $-f''$ would be
very large; recall that $\lim_{\tau\ra 0}(-f'')=\un$.

In view of \re{0.8}, $-f''$ is also dominating $\abs{f'}$, hence $\tilde v$ is a priori
bounded independent of $T$.
\ep

Before we can show that $\norm{Du}$ decays exponentially, we need the following
lemma

\bl\lal{4.4}
For any $k\in\N$ there exists $c_k$ such that
\begin{equation}\lae{4.19}
\nnorm{\h_{\al\bet}}\le c_k \abs\tau ^k.
\end{equation}
Corresponding estimates also hold for $\nnorm{\h_{\al\bet\ga}}, \nnorm{D\psi},
\nnorm{\bar R_{\al\bet}\h^\al}$, or more generally, for any tensor that would vanish
identically, if it would have been formed with respect to the product metric
\begin{equation}
-(dx^0)^2+\bar\s_{ij}dx^idx^j.
\end{equation}
\el

\bp
We only prove the estimate \re{4.19} in detail. The remaining claims can easily be
deduced with the help of the arguments that will follow; in case of $\nnorm{D\psi}$
we use in addition the assumption that all derivatives of $\psi$ of arbitrary order
vanish if $\tau$ tends to $0$.

Let $(\x^\al), (\chi^\al)$ be arbitrary smooth contravariant vectorfields and set
\begin{equation}
\f=\h_{\al\bet}\x^\al\chi^\bet.
\end{equation}

Let us evaluate $\f$ in $(x^0,x)$, $x\in\so$ fixed. Then we have
\begin{equation}
\begin{aligned}
\pde \f{x^0}
=\h_{\al\bet\ga}\x^\al\chi^\bet\h^\ga+\h_{\al\bet}\x^\al_{\hp{\al};\ga}\h^\ga
\chi^\bet +\h_{\al\bet}\x^\al \chi^\bet_{\hp{\bet};\ga}\h^\ga.
\end{aligned}
\end{equation}

\cvm
Since $(\h_{\al\bet})$ is a tensor that vanishes identically in the product metric, we
conclude that $\pde\f{x^0}$ vanishes identically in the product metric, and by
induction we further deduce
\begin{equation}
\lim_{x^0\ra 0}D^k_{x^0}\f =0\qq\A\,k\in\N
\end{equation}
and
\begin{equation}
\abs{D^k_{x^0}\f}\le c_k\qq\A\,k\in \N.
\end{equation}

The mean value theorem then yields
\begin{equation}
\abs{\f(\tau,x)-\f(\tau_0,x)}\le
\sup_{[\tau,\tau_0]}\abs{D_{x^0}\f}\abs{\tau-\tau_0},
\end{equation}
and, by letting $\tau_0$ tend to $0$, we conclude
\begin{equation}
\abs{\f(\tau,x)}\le \sup_{[\tau,0)} \abs{D_{x^0}\f}\msp\abs\tau.
\end{equation}
Applying now induction to $\abs{D_{x^0}\f}$ yields the result because of the
arbitrariness of $(\x^\al), (\chi^\al)$.
\ep

\bl
There exists $\e>0$ and a constant $c_\e$ such that
\begin{equation}
\norm{Du}e^{\e t}\le c_\e\qq\A\,t\in\R[]_+.
\end{equation}
\el

\bp
We employ the relation
\begin{equation}\lae{4.28}
\tilde v^2=1+\norm{Du}^2
\end{equation}
and the fact that $\tilde v$ is uniformly bounded to conclude that for small
$\norm{Du}$
\begin{equation}
2\log \tilde v\sim \norm{Du}^2,
\end{equation}
i.e., we can equivalently prove that $\log \tilde v\msp e^{2\e t}$ is uniformly
bounded.

\cvm
Let $\e>0$ be small and set
\begin{equation}
\f=\log\tilde v\msp e^{2\e t},
\end{equation}
then $\f$ satisfies
\begin{equation}\lae{4.31}
\dot\f -F^{-2}\D\f=\tilde v^{-1}(\dot{\tilde v} -F^{-2}\D\tilde v) e^{2\e t} +
F^{-2}\norm{D\f}^2+ 2\e  \f.
\end{equation}

To get an a priori estimate for $\f$ we shall proceed as in the proof of \rl{4.3}.
For large $T$, $0<T<\un$, assume that
\begin{equation}
\sup_{[0,T]}\sup_{M(t)}\f=\f(t_0,x_0).
\end{equation}
Applying the maximum principle we infer from \re{4.31}, \re{4.2},  \rl{4.2}, and
\rl{4.3}, after multiplying by $F^2$,
\begin{equation}\lae{4.33}
\begin{aligned}
0\le& -\norm A^2e^{2\e t} + c\norm A \abs u e^{2\e t} + c \abs u^2 e^{2\e
t}-nf''\norm {Du}^2 e^{2\e t}\tilde v\\[\cma]
&+c \abs u \norm{Du} e^{2\e t}+c\norm A\msp \norm {Du} e^{2\e t}+ c
\norm{Du}^2 e^{2\e t} +2 \e F^2\f.
\end{aligned}
\end{equation}

Now, we have
\begin{equation}\lae{4.34}
\begin{aligned}
F^2&=H^2+n^2\abs{f'}^2\tilde v^2+n^2\abs{\psi_\al\nu^\al}^2\\[\cma]
&\hp{=}-2nHf'\tilde v+2nH\psi_\al\nu^\al-2n^2 f' \tilde v \psi_\al\nu^\al,
\end{aligned}
\end{equation}
hence $\f$ is apriori bounded, if $\e$ is small enough, $0<\e<<\tilde\ga$.

\cvm
Here we also used the boundedness of $\tilde v$ so that
\begin{equation}
\f\le c e^{2\e t},
\end{equation}
as well as the boundedness of $\tilde u=u e^{\ga t}$.

To control the term
\begin{equation}
\e n^2 \abs{f'}^2\tilde v^2 \f
\end{equation}
we employed the assumption \re{0.8} yielding
\begin{equation}\lae{4.37}
-c\le f''+\tilde \ga \abs{f'}^2\le c
\end{equation}
as well as the estimate
\begin{equation}\lae{4.38}
\abs{\log\tilde v-\tfrac12\norm{Du}^2}\le c \norm{Du}^4
\end{equation}
because of \re{4.28}.
\ep

After having established the exponential decay of $\norm{Du}$, we can improve the
decay rate.

\bl\lal{4.6}
For any $0<\lam<\ga$ there exists $c_\lam$ such that
\begin{equation}
\norm{Du} e^{\lam t}\le c_\lam.
\end{equation}
\el

\bp
As in the proof of the preceding lemma set
\begin{equation}
\f=\log\tilde v\msp[2] e^{2\lam t}.
\end{equation}
Let $T$, $0<T<\un$, be large and $(t_0,x_0)$ be such that
\begin{equation}
\sup_{[0,T]}\sup_{M(t)}\f=\f(t_0,x_0).
\end{equation}
Applying the maximum principle we then obtain an inequality as in \re{4.33}, where
$\e$ has to be replaced by $\lam$.

\cvm
The bad terms which need further consideration are part of
\begin{equation}
2\lam F^2\f,
\end{equation}
especially
\begin{equation}\lae{4.43}
2\lam H^2\f
\end{equation}
and 
\begin{equation}
2\lam n^2 \abs{f'}^2\tilde v^2\f.
\end{equation}

The quantity in \re{4.43} can be absorbed by
\begin{equation}
-\norm A^2 e^{2\lam t},
\end{equation}
since $\f=\log\tilde v\msp[2] e^{2\lam t}$ and $\log\tilde v$ decays exponentially.

The second term is dominated by
\begin{equation}
-nf''\norm{Du}^2 e^{2\lam t}\tilde v,
\end{equation}
because of \re{4.28}, \re{4.37}, \re{4.38}, the exponential decay of $\norm{Du}$,
and the assumption that $\lam <\ga$.

Thus we see that $\f$ is a priori bounded independent of $T$.
\ep

\section{$C^2$-estimates}\las{5}

The ultimate goal is to show that $\norm A e^{\ga t}$ is uniformly bounded.
However, this result can only be derived by first establishing some preliminary
estimates.

Let us start by proving that $F$ grows exponentially fast. From the evolution
equation \re{2.11} we deduce
\begin{equation}\lae{5.1}
\dot H-F^{-2}\D F=-2F^{-3}\norm{DF}^2+F^{-2}(\norm A^2+\bar
R_{\al\bet}\nu^\al\nu^\bet)F,
\end{equation}
where we have used that
\begin{equation}
\dot H=\de^i_j\dot h^j_i.
\end{equation}

Replacing $\dot H$ by $\dot F$ in the evolution equation \re{5.1} and observing that
\begin{equation}
\begin{aligned}
\dot F&= \dot H -n f''\tilde v^2F^{-1}+nf'\h_{\al\bet}\nu^\al\nu^\bet
F^{-1}\\[\cma]
&\hp{=}+nf'u^iF_iF^{-2}-n\psi_{\al\bet}\nu^\al\nu^\bet F^{-1}+n\psi_\al
x^\al_iF^iF^{-2}
\end{aligned}
\end{equation}
we obtain
\begin{equation}\lae{5.4}
\begin{aligned}
\dot F&-F^{-2}\D F=-2 F^{-3}\norm{DF}^2 +F^{-2}(\norm A^2 +\bar
R_{\al\bet}\nu^\al\nu^\bet)F\\[\cma]
&+F^{-2}(-nf''\tilde v^2+nf'\h_{\al\bet}\nu^\al\nu^\bet-n\psi_{\al\bet}
\nu^\al\nu^\bet)F\\[\cma]
&+F^{-2}(nf'u_i+n\psi_\al x^\al_i)F^i.
\end{aligned}
\end{equation}

\bl
There exist positive constants $\de$ and $c_\de$ such that
\begin{equation}
c_\de e^{\de t}\le F\qq\A\,t\in \R[]_+.
\end{equation}
\el

\bp
Define
\begin{equation}
\f=Fe^{-\de t}.
\end{equation}
Let $T$, $0<T<\un$, be large and $(t_0,x_0)$ be such that
\begin{equation}
\sup_{[0,T]}\sup_{M(t)}\f=\f(t_0,x_0).
\end{equation}
Applying the maximum principle we deduce from \re{5.4}
\begin{equation}
\begin{aligned}
0&\ge \norm A^2 +\bar R_{\al\bet}\nu^\al\nu^\bet + nf'\h_{\al\bet}
\nu^\al\nu^\bet-n f''\tilde v-n\psi_{\al\bet}\nu^\al \nu^\bet -\de F^2,
\end{aligned}
\end{equation}
and we further conclude that, for small $\de$, $t_0$ cannot exceed a certain value in
view of the relations \re{4.34} and \re{4.37}, hence the result.
\ep

Replacing in \re{5.1} $F$ by $H$ we obtain an evolution equation for $H$
\begin{equation}\lae{5.9}
\begin{aligned}
\dot H&-F^{-2}\D H=-2F^{-3}\norm{DF}^2+F^{-2} (\norm A^2+\bar
R_{\al\bet}\nu^\al \nu^\bet)F\\[\cma]
&+F^{-2}(nf''\tilde v^2 H-nf''\tilde v g^{ij}\bar h_{ij}-nf'''\norm{Du}^2\tilde v - 2n
f''\h_{\al\bet} \nu^\al x^\bet_i u^i\\[\cma]
&\q+4nf'' h_{ij}u^iu^j -nf'\h_{\al\bet\ga}\nu^\al x^\bet_i x^\ga_jg^{ij} -2nf'h^{ij}
\h_{\al\bet} x^\al_ix^\bet_j\\[\cma]
&\q -n\h_{\al\bet}\nu^\al\nu^\bet f' H -n f'\norm A^2 \tilde v +nf' H_ku^k +
nf'\bar R_{\al\bet} \nu^\al x^\bet_ku^k)\\[\cma]
&+nF^{-2}(\psi_{\al\bet\ga} \nu^\al x^\bet_ix^\ga_jg^{ij}+\psi_{\al\bet}
\nu^\al\nu^\bet H+2\psi_{\al\bet} x^\al_ix^\bet_j h^{ij}\\[\cma]
&\q+\norm A^2\psi_\al\nu^\al+\psi_\al x^\al_k H^k+\bar R_{\al\bet}\nu^\al
x^\bet_k\psi_\ga x^\ga_l  g^{kl}).
\end{aligned}
\end{equation}

In deriving this equation we used the Weingarten and Codazzi equations, the definition
of $F$ and the relation
\begin{equation}
\tilde vH=-\D u+g^{ij}\bar h_{ij},
\end{equation}
where $\bar h_{ij}$ is the second fundamental form of the slices $\{x^0=\const\}$.

\bl
$H$ is uniformly bounded from below during the evolution.
\el

\bp
Let $T$, $0<T<\un$, be large and $x_0=x(t_0,\x_0)$ be such that
\begin{equation}
\inf_{[0,T]}\inf_{M(t)}H=H(x_0).
\end{equation}
Applying the maximum principle and some trivial estimates we deduce from \re{5.9}
\begin{equation}
\begin{aligned}
0&\ge -2F^{-3}\norm{DF}^2+F^{-2}(\norm A^2+\bar R_{\al\bet}\nu^\al
\nu^\bet)F\\[\cma]
&\hp{\ge} +F^{-2}(nf''\tilde vH-c\abs{f'}^\frac32 -\tfrac n2f'\norm A^2\tilde
v-c(1+\norm A^2)),
\end{aligned}
\end{equation}
where we have used \rc{3.7}, \rl{4.4}, \rl{4.6} and assumed that $H(x_0)\le -1$.

\cvm
To estimate the term involving $\norm{DF}^2$ we note that
\begin{equation}
\begin{aligned}
\norm{DF}^2&= \norm{DH}^2+n^2\abs{f''}^2\norm{Du}^2\tilde
v^2+n^2\abs{f'}^2\norm{D\tilde v}^2\\[\cma]
&\hp{=}+n^2\norm{D(\psi_\al\nu^\al)}^2-2nf''H_ku^k\tilde
v-2nf'H_k\tilde v^k\\[\cma]
&\hp{=}+2nH^k(\psi_\al\nu^\al)_k + 2n^2f'f''\tilde v_ku^k\tilde v\\[\cma]
&\hp{=}-2n^2f''\tilde v(\psi_\al\nu^\al)_ku^k-2n^2 f'(\psi_\al\nu^\al)_k\tilde v^k.
\end{aligned}
\end{equation}

$DH$ vanishes in $x_0$, and because of \re{4.4}, \rl{4.4} and \rl{4.6} we have
\begin{equation}
\begin{aligned}
\norm{D\tilde v}\le c\nnorm{\h_{\al\bet}}+\norm A\norm {Du}\le c_\lam (1+\norm
A)e^{-\lam t}\q\A\,0<\lam<\ga.
\end{aligned}
\end{equation}

\cvm
Combining these estimates with the exponential growth of $F$ we conclude
\begin{equation}
F^{-1}\norm{DF}^2\le c(1+\abs{f''}+\norm A^2),
\end{equation}
hence the a priori bound from below for $H$.
\ep

Next we shall show that the principal curvatures of $M(t)$ are uniformly bounded
from above, i.e., we want to estimate $h^j_i$ from above.

\cvm
Let us first derive a parabolic equation satisfied by $h^j_i$ from the evolution
equation \re{2.11}.

\cvm
Using the definition of $F$ we immediately obtain
\begin{equation}
\begin{aligned}
\dot h^j_i-F^{-2}
H^j_i&=-2F^{-3}F_iF^j+F^{-1}h_{ik}h^{kj}\\[\cma]
&\hp{=}+F^{-1}\riema\al\bet\ga\de \nu^\al
x^\bet_i\nu^\ga x^\de_k g^{kj}\\[\cma]
&\hp{=}+F^{-2}(-n(f'\tilde v)^j_i+n(\psi_\al\nu^\al)^j_i)
\end{aligned}
\end{equation}
and conclude further
\begin{equation}\lae{5.17}
\begin{aligned}
\dot h^j_i&-F^{-2}\D
h^j_i=\\[\cma]
&-2F^{-3}F_iF^j+F^{-1}h_{ik}h^{kj}+F^{-1}\riema\al\bet\ga\de \nu^\al
x^\bet_i\nu^\ga x^\de_k g^{kj}\\[\cma]
&-F^{-2}\norm A^2h^j_i+F^{-2}Hh_{ik}h^{kj}+2F^{-2}h^{kl}\riema \al\bet\ga\de
x^\al_kx^\bet_ix^\ga_lx^\de_rg^{rj}\\[\cma]
&-F^{-2}(g^{kl}\riema\al\bet\ga\de
x^\al_mx^\bet_kx^\ga_rx^\de_lh^m_ig^{rj}+g^{kl}\riema\al\bet\ga\de
x^\al_mx^\bet_kx^\ga_ix^\de_lh^{mj}\\[\cma]
&\q +\bar R_{\al\bet}\nu^\al\nu^\bet h^j_i-H\riema\al\bet\ga\de \nu^\al
x^\bet_i\nu^\ga x^\de_mg^{mj})\\[\cma]
&+F^{-2}g^{kl}\bar R_{\al\bet\ga\de;\e}(\nu^\al
x^\bet_kx^\ga_lx^\de_ix^\e_mg^{mj}+\nu^\al x^\bet_ix^\ga_kx^\de_mx^\e_l
g^{mj})\\[\cma]
&+F^{-2}(nf''h^j_i\tilde v^2+nf''\tilde v \h_{\al\bet}x^\al_ix^\bet_kg^{kj}
-nf'''u_iu^j\tilde v\\[\cma]
&\q -nf''(\tilde v_iu^j+\tilde v^ju_i)-nf'[\h_{\al\bet\ga}\nu^\al
x^\bet_ix^\ga_kg^{kj}+\h_{\al\bet}x^\al_kx^\bet_ih^{kj}\\[\cma]
&\q+\h_{\al\bet}\nu^\al\nu^\bet h^j_i+h^k_ih^j_k\tilde
v-h^j_{i;k}u^k+\riema\al\bet\ga\de \nu^\al x^\bet_ix^\ga_kx^\de_l
g^{lj}u_k])\\[\cma]
&+nF^{-2}(\psi_{\al\bet\ga}\nu^\al x^\bet_i x^\ga_kg^{kj}+
\psi_{\al\bet}\nu^\al\nu^\bet h^j_i +\psi_{\al\bet} x^\al_k x^\bet_i
h^{kj}\\[\cma] 
&\q+\psi_{\al\bet}x^\bet_lx^\al_k h^k_ig^{lj}+\psi_\al\nu^\al
h_{ki}h^{kj}+\psi_\al x^\al_k h^k_{i;l}g^{lj}),
\end{aligned}
\end{equation}
where we used the relation
\begin{equation}\lae{5.18}
h_{ij}\tilde v=-u_{ij}+\bar h_{ij}=-u_{ij}-\h_{\al\bet}x^\al_i x^\bet_j,
\end{equation}
equation \re{4.4} as well as the Weingarten and Codazzi equations.

\bl
The principal curvatures $\ka_i$ of $M(t)$ are uniformly bounded during the evolution.
\el

\bp
Since we already know that $H\ge -c$, it suffices to prove an uniform estimate from
above.

Let $\f$ be defined by
\begin{equation}
\f=\sup\set{h_{ij}\h^i\h^j}{\norm\h=1}.
\end{equation}
We shall prove that
\begin{equation}
w=\log\f+\lam \tilde v
\end{equation}
is uniformly bounded from above, if $\lam$ is large enough.

Let $0<T<\un$ be large, and $x_0=x_0(t_0)$, with $ 0<t_0\le T$, be a point in
$M(t_0)$  such that
\begin{equation}
\sup_{M_0}\f<\sup\set {\sup_{M(t)} \f}{0<t\le T}=\f(x_0).
\end{equation}

We then introduce a Riemannian normal coordinate system $(\x^i)$ at $x_0\in
M(t_0)$ such that at $x_0=x(t_0,\x_0)$ we have
\begin{equation}
g_{ij}=\de_{ij}\q \tup{and}\q \f=h_n^n.
\end{equation}

Let $\tilde \h=(\tilde \h^i)$ be the contravariant vector field defined by
\begin{equation}
\tilde \h=(0,\dotsc,0,1),
\end{equation}
and set
\begin{equation}
\tilde \f=\frac{h_{ij}\tilde \h^i\tilde \h^j}{g_{ij}\tilde \h^i\tilde \h^j}\raise 2pt
\hbox{.}
\end{equation}

$\tilde \f$ is well defined in neighbourhood of $(t_0,\x_0)$, and $\tilde \f$
assumes its maximum at $(t_0,\x_0)$. Moreover, at $(t_0,\x_0)$ we have
\begin{equation}
\dot{\tilde \f}=\dot h_n^n,
\end{equation}
and the spatial derivatives do also coincide; in short, at $(t_0,\x_0)$ $\tilde \f$
satisfies the same differential equation \re{5.17} as $h_n^n$. For the sake of
greater clarity, let us therefore treat $h_n^n$ like a scalar and pretend that
$w$ is defined by
\begin{equation}
w=\log h^n_n+\lam \tilde v.
\end{equation}

At $(t_0,\x_0)$ we have $\dot w\ge 0$, and, in view of the maximum principle, we
deduce from \re{5.17} and \re{4.2}
\begin{equation}
\begin{aligned}
0&\le -\lam\norm A^2\tilde v+c\lam (1+\norm A+\abs{f'}\msp e^{-\frac32 \ga
t}\norm A)\\[\cma]
&\hp{\le}+c(\abs{H}h^n_n+\abs{f'}h^n_n+\abs{f'})+nf''\tilde v^2\\[\cma]
&\hp{\le}+c\abs{f'}\msp \norm{D\log h^n_n}\msp[2] \norm{Du}+\norm{D\log
h^n_n}^2+c\norm{D\log h^n_n},
\end{aligned}
\end{equation}
where we assumed  $h^n_n\ge 1$, and in addition used \re{4.4} and the known
exponential decay estimates for $\norm{Du}$.

\cvm
Since $Dw=0$ in $x_0$, we have
\begin{equation}
\norm{D\log h^n_n}=\lam\norm{D\tilde v}\le \lam c (1+\norm A\msp[1]
\norm{Du}).
\end{equation}
Hence, if $\lam$ is chosen large enough, we obtain an a priori bound for $h^n_n$
from above.
\ep

An immediate corollary is

\bc\lac{5.4}
There exist positive constants $c_1,c_2$ such that
\begin{equation}
c_1\le Fe^{-\ga t}\le c_2.
\end{equation}
\ec

\bp
Since $H$ is uniformly bounded we conclude
\begin{equation}\lae{5.30}
Fe^{-\ga t}\sim -nf'e^{-\ga t}\tilde v=-nf'u(ue^{\ga t})^{-1}\tilde v
\end{equation}
and the result follows from \rl{3.4} and \rt{3.6}.
\ep

We can now prove an exponential decay for $\norm A$.

\bl
For any $0<\lam <\ga$ there exists $c_\lam$ such that
\begin{equation}
\norm A e^{\lam t}\le c_\lam\qq\A\,t\in\R[]_+.
\end{equation}
\el

\bp
Let $\f=\frac12\norm A^2$, then
\begin{equation}\lae{5.32}
\dot\f-F^{-2}\D\f=-F^{-2}\norm{DA}^2+(\dot h^j_i-F^{-2}\D h^j_i)h^i_j,
\end{equation}
where
\begin{equation}
\norm{DA}^2=h_{ij;k}h^{ij}_{\hp{ij};l}g^{kl}.
\end{equation}

\cvm
Define $w=\f e^{2\lam t}$ with $0<\lam<\ga$. Let $0<T<\un$ be large, and 
$x_0=x_0(t_0)$, with $ 0<t_0\le T$, be a point in
$M(t_0)$  such that
\begin{equation}
\sup_{M_0}w<\sup\set {\sup_{M(t)} w}{0<t\le T}=w(x_0).
\end{equation}

Applying the maximum principle we deduce from \re{5.32} and \re{5.17}
\begin{equation}\lae{5.35}
\begin{aligned}
0&\le -\norm{DA}^2e^{2\lam t} -F^{-1}h^{ij}F_iF_j e^{2\lam t}+2nf''\tilde
v^2w\\[\cma]
&\hp{\le}+c_\e e^{-\e t}\abs{f''} \msp\norm A e^{\lam t} +c\abs{f'}(w+1)+2\lam
F^2 w,
\end{aligned}
\end{equation}
with some small positive $\e=\e(\lam)$; here we used \rl{4.6} and \rc{3.7}.

\cvm
It remains to estimate the second and the last term in the preceding inequality. The
only relevant term in $2\lam F^2 w$ is
\begin{equation}
2\lam n^2\abs{f'}^2\tilde v^2 w;
\end{equation}
combining it with $2nf''\tilde v^2 w$ gives
\begin{equation}
2nf''\tilde v^2 w+2\lam n^2\abs{f'}^2\tilde v^2 w\le -2n^2(\ga
-\lam)\abs{f'}^2\tilde v^2w+cw,
\end{equation}
in view of \re{4.37}.

The remaining term can be estimated
\begin{equation}
-F^{-1}h^{ij}F_iF_je^{2\lam t}\le c e^{-\ga t} \norm {DA}^2 e^{2\lam
t}+c_\e\abs{f'}^2e^{-\e t}\norm A e^{\lam t} +c(1+w),
\end{equation}
with some positive $\e=\e(\lam)$.

Inserting these estimates in \re{5.35} we obtain an a priori bound for $w$.
\ep

Though we now could prove an a priori estimate for $\norm A e^{\ga t}$, let us first
derive a corresponding estimate for $\norm{Du} e^{\ga t}$. The estimate for the
second fundamental form is then slightly easier to prove.

\bt\lat{5.6}
Let $\tilde u=u^{\ga t}$, then $\norm{D\tilde u}$ is uniformly bounded during the
evolution.
\et

\bp
Let $\f=\f(t)$ be defined by
\begin{equation}
\f=\sup_{M(t)}\log\tilde v\msp[2]  e^{2\ga t}.
\end{equation}

Then, in view of the maximum principle, we deduce from equation \re{4.2}
\begin{equation}\lae{5.40}
\dot\f\le c e^{-\e t}+F^{-2}(nf''\norm{D\tilde u}^2\tilde v+2\ga F^2 w)
\end{equation}
for some positive $\e$, where we haved used the known exponential decay of
$\norm A$ and $\norm{Du}$ as well as \rl{4.2}, \rl{4.4}, \rc{5.4} and the
inequalities \re{4.37} and \re{4.38}; the inequality is valid for \aev $t$.

The second term on the right-hand side of \re{5.40} can be estimated from above by
\begin{equation}
c e^{-\e t}(1+w),
\end{equation}
in view of \re{4.37}, \re{4.38} and the known decay of $\norm A$, $\norm{Du}$ as
well as the result in \rc{5.4}. Hence we conclude
\begin{equation}
\dot \f\le c e^{-\e t}(1+\f),
\end{equation}
i.e., $\f$ is uniformly bounded.
\ep

\bt\lat{5.7}
The quantity $w=\frac12 \norm A^2 e^{2\ga t}$ is uniformly bounded during the
evolution.
\et

\bp
Define $\f=\f(t)$ by
\begin{equation}
\f=\sup_{M(t)} w.
\end{equation}

Applying the maximum principle we deduce from \re{5.32} that for \aev $t$
\begin{equation}
\begin{aligned}
\dot\f&\le -F^{-2}\norm{DA}^2 e^{2\ga t}+F^{-3}(-2 h^{ij}F_iF_je^{2\ga t}-nf'''
h^{ij}u_iu_j\tilde v)\\[\cma]
&\hp{\le}+F^{-2}(n f'' \tilde v^2 \f +\ga F^2\f)+ce^{-\e t}(1+\f)\\[\cma]
&\hp{\le}+F^{-1}\riema
\al\bet\ga\de \nu^\al x^\bet_i\nu^\ga x^\de_j h^{ij} e^{2\ga t}
\end{aligned}
\end{equation}

\cvm
The last terms on the right-hand side of this inequality can be estimated as follows
\begin{equation}\lae{5.45}
\begin{aligned}
F^{-3}&(-2 h^{ij}F_iF_je^{2\ga t}-nf'''
h^{ij}u_iu_j\tilde v)\le\\[\cma]
&F^{-3}(-2 \abs{f''}^2+f'f''')h^{ij}\tilde u_i\tilde u_j\tilde v^2
n^2+cF^{-3}\norm{DA}^2e^{2\ga t}\\[\cma]
&+ce^{-\e t}(1+\f).
\end{aligned}
\end{equation}

Now, we observe that
\begin{equation}
(f''+\tilde \ga \abs{f'}^2)'=f'''+2\tilde \ga f'f''={\mc C} f',
\end{equation}
where $\mc C$ is a bounded function in view of assumption \re{0.9}, and hence
\begin{equation}
2\abs{f''}^2-f'f'''=2\abs{f''}^2+2\tilde\ga \abs{f'}^2f''-\mc C\abs{f'}^2,
\end{equation}
i.e.,
\begin{equation}
\abs{2\abs{f''}^2-f'f'''}\le c \abs{f'}^2,
\end{equation}
and we conclude that the left-hand side of \re{5.45} can be estimated from above by
\begin{equation}
ce^{-\e t} (1+\f) +cF^{-2}\norm{DA}^2 e^{\ga t}
\end{equation}

\cvm
Next, we estimate
\begin{equation}
F^{-2}(nf''\tilde v^2 +\ga F^2)\f\le ce^{\-\e t}\f,
\end{equation}
and finally
\begin{equation}
\begin{aligned}
F^{-1}\riema
\al\bet\ga\de \nu^\al x^\bet_i\nu^\ga x^\de_j h^{ij} e^{2\ga t}&\le ce^{-\e
t}(1+\f) +F^{-1}\riema 0i0j h^{ij}e^{2\ga t} \tilde v^2,
\end{aligned}
\end{equation}
but
\begin{equation}
\abs{\riema 0i0j}\le c \abs{u},
\end{equation}
\cf \rl{4.4}.

Hence, we deduce
\begin{equation}
\dot \f \le ce^{-\e t}(1+\f)
\end{equation}
for some positive $\e$ and for \aev $t$, i.e., $\f$ is bounded.
\ep

\section{Higher order estimates}\las{6}

After having established the boundedness of
\begin{equation}
\norm A^2 e^{2\ga t}
\end{equation}
corresponding estimates for the derivatives of the second fundamental form will be
proved recursively.

Our starting point is the equation \re{5.17}. It contains two very bad terms
\begin{equation}
-nF^{-2}f'''u_iu^j\tilde v,
\end{equation}
and another one which is hidden in the expression
\begin{equation}
-2F^{-3}F_iF^j.
\end{equation}

\cvm
To handle these terms we proceed as in the proof of \rt{5.7} by combining the two
crucial terms in
\begin{equation}\lae{6.4}
F^{-3}(-2F_iF^j-nFf'''u_iu^j\tilde v)
\end{equation}
to
\begin{equation}\lae{6.5}
F^{-3}(-2\abs{f''}^2+f'f''')u_iu^jn^2\tilde v^2
\end{equation}
and observing that
\begin{equation}
\f=-2\abs{f''}^2+f'f'''=(f''+\tilde\ga\abs{f'}^2)'f'-2f''(f''+\tilde\ga\abs{f'}^2).
\end{equation}

In view of our assumption \re{0.9} and \rc{3.7} we conclude that the spatial
derivatives of
$\f$ can be estimated by
\begin{equation}\lae{6.7}
\norm{D^m\f}\le c_m\norm{\tilde u}_m e^{2\ga t}\qq\A\,m\in\N.
\end{equation}

\cvm
Let us introduce the following abbreviations

\bd
(i) For arbitrary tensors $S,T$ denote by $S\star T$ any linear combination of
tensors formed by contracting over $S$ and $T$. The result can be a tensor or a
function. Note that we do not distinguish between $S\star T$ and $c\msp S\star T$,
$c$ a constant.

\cvm
(ii) The symbol $A$ represents the second fundamental form of the hypersurfaces
$M(t)$ in $N$, $\tilde A=A e^{\ga t}$ is the scaled version, and $D^mA$ \resp
$D^m\tilde A$ represent the covariant derivatives of order $m$.

\cvm
(iii) For $m\in \N$ denote by $\mc O_m$ a tensor expression defined on $M(t)$ that
satisfies the pointwise estimates
\begin{equation}
\norm{\mc O_m}\le c_m(1+\norm{\tilde A}_m)^{p_m},
\end{equation}
where $c_m,p_m$ are positive constants, and
\begin{equation}
\norm{\tilde A}_m=\sum_{\abs\al \le m}\norm{D^\al \tilde A}.
\end{equation}
Moreover, the derivative of $\mc O_m$ is of class $\mc O_{m+1}$ and can be
estimated by
\begin{equation}
\norm{D\mc O_m}\le c_m(1+\norm{\tilde A}_m)^{p_m}(1+\norm{D^{m+1}\tilde
A})
\end{equation}
with (different) constants $c_m,p_m$.

\cvm
(iv) The symbol $\mc O$ represents a tensor such that $D\mc O$ is of class $\mc
O_0$.
\ed

\br
We emphasize the following relations
\begin{equation}
D^m\mc O_0=\mc O_m\qq\A\,m\in\N,
\end{equation}
\begin{equation}
F^{-1}DF=F^{-1}DA+\mc O,
\end{equation}
\begin{equation}
DF e^{-\ga t}=e^{-\ga t} DA +\mc O,
\end{equation}
\begin{equation}
F^{-1}\mc O_m=\mc O_m\qq\A\,m\in\N,
\end{equation}
and
\begin{equation}\lae{6.15}
\abs{\riema 0i0j}\le c_m\abs{u}^m\qq\A\,m\in\N,
\end{equation}
\cf \rl{4.4}.
\er

With these definitions and the relations \re{6.5} and \re{6.7} in mind we can write the
evolution equation for $\tilde h^j_i$ in the form
\begin{equation}\lae{6.16}
\begin{aligned}
\dot{\tilde h}^j_i-F^{-2}\D \tilde h^j_i&=F^{-3}D\tilde A\star DA+F^{-2}\mc
O\star D\tilde A\\[\cma]
&\hp{=}+F^{-3}\mc O_0\star D\tilde A +F^{-2}\mc O_0+F^{-1}\mc O,
\end{aligned}
\end{equation}
where the right-hand side is considered to be a mixed tensor of order two though we
omitted the indices.

\cvm
Using the fact that
\begin{equation}\lae{6.17}
\dot g_{ij}=-2F^{-1}h_{ij}=-2F^{-1}e^{-\ga t}\tilde h_{ij}=F^{-2}\mc O_0
\end{equation}
we can rewrite \re{6.16} in the form
\begin{equation}\lae{6.18}
\begin{aligned}
\dot{\tilde A}-F^{-2}\D \tilde A&=F^{-3}D\tilde A\star DA+F^{-2}\mc
O\star D\tilde A\\[\cma]
&\hp{=}+F^{-3}\mc O_0\star D\tilde A +F^{-2}\mc O_0+F^{-1}\mc O
\end{aligned}
\end{equation}
regardless of  representing  $\tilde A$ as a covariant, contravariant or mixed
tensor.

\cvm
Differentiating this equation covariantly with respect to a spatial variable we deduce
\begin{equation}\lae{6.19}
\begin{aligned}
\tfrac D{dt}(D\tilde A)&-F^{-2}\D D\tilde A=F^{-1}\mc O_0+F^{-3}D^2\tilde A\star
DA+F^{-2}\mc O\star D^2\tilde A\\[\cma]
&+F^{-4}D\tilde A\star DA\star DA+F^{-3}\mc O\star D\tilde A\star
DA+F^{-2}D\tilde A\star \mc O_0\\[\cma]
&+F^{-4}D\tilde A\star DA\star \mc O_0+F^{-3}D\tilde A\star D\mc
O_0+F^{-3}D^2\tilde A\star \mc O_0,
\end{aligned}
\end{equation}
where we  used the Ricci identities to commute the second derivatives of a tensor.

\cvm
Finally, using induction, we conclude
\begin{equation}\lae{6.20}
\begin{aligned}
\tfrac D{dt}(D^{m+1}\tilde A)&-F^{-2}\D D^{m+1}\tilde A=F^{-1}\mc
O_m+F^{-3}D^{m+2}\tilde A\star DA\\[\cma]
&+F^{-2}D^{m+1}A\star \mc O_m+F^{-3}D^{m+2}\tilde A\star \mc O_0\\[\cma]
&+\Theta F^{-3}D^{m+1}\tilde A\star D^{m+1}A,
\end{aligned}
\end{equation}
for any $m\in\N^*$, where $\Theta =1$, if $m=1$, and $\Theta=0$ otherwise.

\cvm
We are now going to prove uniform bounds for $\tfrac12 \norm{D^{m+1}\tilde
A}^2$ for all $m\in\N$.

First we observe that
\begin{equation}
\begin{aligned}
\tfrac D{dt}&(\tfrac12 \norm{D^{m+1}\tilde
A}^2)-F^{-2}\D \tfrac12 \norm{D^{m+1}\tilde
A}^2=
-F^{-2}\norm{D^{m+2}\tilde A}^2\\[\cma]
&+F^{-1}\mc
O_m\star D^{m+1}\tilde A +F^{-3}D^{m+2}\tilde A\star DA\star D^{m+1}\tilde
A\\[\cma] 
&+F^{-2}D^{m+1}A\star \mc O_m\star D^{m+1}\tilde A+F^{-3}D^{m+2}\tilde
A\star \mc O_0\star D^{m+1}\tilde A\\[\cma] 
&+\Theta F^{-3}D^{m+1}\tilde A\star D^{m+1}A\star D^{m+1}\tilde A,
\end{aligned}
\end{equation}
if $m\in\N^*$, in view of \re{6.20}, where similar equations are also valid for
$\frac12\norm{\tilde A}^2$ and $\frac12\norm {D\tilde A}^2$, \cf \re{6.18} and
\re{6.19}.

\bt\lat{6.3}
The quantities $\frac12\norm{D^m\tilde A}^2$ are uniformly bounded during the
evolution for all $m\in\N^*$.
\et

\bp
We proof the theorem recursively by estimating
\begin{equation}
\f=\log \tfrac12\norm{D^{m+1}\tilde A}^2+\m\tfrac12\norm{D^m\tilde
A}^2+\lam e^{-\ga t},
\end{equation}
where $\m$ is a small positive constant
\begin{equation}
0<\m=\m(m)<<1,
\end{equation}
and $\lam$ large, $\lam=\lam(m)>>1$.

We shall only treat the case $m=0$, since then the structure of the right-hand side is
worst, at least formally, \cf \re{6.19}.

\cvm
Fix $0<T<\un$, $T$ very large, and suppose that
\begin{equation}
2\sup\norm{\tilde A}^2<\sup_{[0,T]}\sup_{M(t)}\f=\f(x(t_0,\x_0))
\end{equation}
for $0<t_0\le T$, where $e^{-\ga t_0}$ should be small compared with $\m$, i.e.,
$t_0$ has to be large.

Applying the maximum principle we deduce
\begin{equation}\lae{6.25}
\begin{aligned}
0&\le \m^2 F^{-2}\norm{D\tfrac12\norm{\tilde A}^2}-F^{-2}\norm{D^2\tilde
A}^2\norm{\tilde A}^{-2}-\tfrac\lam2\ga e^{-\ga t}\\[\cma]
&\hp{\le} -\tfrac\m2 F^{-2}\norm{D\tilde A}^2+c F^{-4}\norm{D\tilde A}^2.
\end{aligned}
\end{equation}

Now, we observe that
\begin{equation}
\norm{D\tfrac12\norm{\tilde A}^2}\le c\norm{D\tilde A}^2\msp \norm{\tilde
A}^2\le c\norm{D\tilde A}^2
\end{equation}
and hence the right-hand side of inequality \re{6.25} would be negative, if $\m$ is
small, $\lam$ large and $t_0$ large.

Thus $\f$ is a priori bounded.

The proof for $m\ge 1$ is similar.
\ep

\section{Convergence of $\tilde u$ and the behaviour of derivatives in $t$}\las{7}

Let us first prove that $\tilde u$ converges when $t$ tends to infinity.

\bl\lal{7.1}
$\tilde u$ converges in $C^m(\so)$ for any $m\in\N$, if $t$ tends to infinity, and
hence $D^m\tilde A$ converges.
\el

\bp
$\tilde u$ satisfies the evolution equation
\begin{equation}
\begin{aligned}
\dot{\tilde u}=\frac{\tilde ve^{\ga t}}F +\ga \tilde u=\frac{\tilde ve^{\ga
t}}F(1-\tilde\ga f'u+v\ga H e^{-\ga t}+v\ga n\psi_\al\nu^\al e^{-\ga t}),
\end{aligned}
\end{equation}
hence we deduce
\begin{equation}
\abs{\dot{\tilde u}}\le c e^{-2\ga t},
\end{equation}
in view of \rl{3.4} and the known estimates for $H,F$ and $\psi$, i.e., $\tilde u$
converges uniformly. Due to \rt{6.3}, $D^m\tilde u$ is uniformly bounded, hence
$\tilde u$ converges in $C^m(\so)$.

\cvm
The convergence of $D^m\tilde A$ follows from \rt{6.3} and  the convergence of
$\tilde h_{ij}$, which in turn can be deduced from equation \re{5.18}.
\ep

Combining the equations \re{6.18}, \re{6.19}, \re{6.20}, and \rt{6.3} we
immediately conclude

\bl\lal{7.2}
$\norm{\frac D{dt}D^m\tilde A}$ and $\norm{\frac D{dt} D^mA}$ decay by the
order
$e^{-\ga t}$ for any $m\in\N$.
\el

\bc\lac{7.3}
$\tfrac D{dt}D^mA e^{\ga t}$ converges, if $t$ tends to infinity.
\ec

\bp
Applying the product rule we obtain
\begin{equation}
\tfrac D{dt}D^m\tilde A=\tfrac D{dt}D^m Ae^{\ga t}+\ga D^m\tilde A,
\end{equation}
hence the result, since the left-hand side converges to zero and $D^m\tilde A$
converges.
\ep

In view of \rl{3.4} $f'u$ converges to $\tilde\ga^{-1}$, if $t$ tends to infinity,
moreover, because of the condition \re{0.10} and the estimates for $u$ \resp $\tilde
u$, we further deduce

\bl\lal{7.4}
For any $m\in\N$ we have
\begin{equation}
\norm{D^m(f'u)}\le c_m.
\end{equation}
\el

\bp
We only consider the case $m=1$. Differentiating $f'u$ we get
\begin{equation}
(f'u)_k=f''uu_k+f'u_k=f''u^2u^{-1}u_k+f'uu^{-1}u_k,
\end{equation}
but
\begin{equation}
u^{-1}u_k=\tilde u^{-1}\tilde u_k
\end{equation}
and hence uniformly bounded in view of \rt{3.6} and \rt{5.6}.
\ep

\bc\lac{7.4}
We have
\begin{equation}\lae{7.6}
\norm{D^mF^{-1}}\le c_m F^{-1}\qq\A\,m\in\N.
\end{equation}
\ec

\bp
Recall that
\begin{equation}
F=H-n\tilde vf'+n\psi_\al\nu^\al
\end{equation}
and hence
\begin{equation}
(F^{-1})_k=-F^{-2}(H_k-n\tilde v_kf'-n\tilde vf''u_k+n(\psi_\al\nu^\al)_k).
\end{equation}

Now, writing
\begin{equation}
\begin{aligned}
F^{-1}&(H_k-n\tilde v_kf'-n\tilde vf''u_k+n(\psi_\al\nu^\al)_k)=\\[\cma]
&(Fu)^{-1}(uH_k-n\tilde v_kf'u-n\tilde vf''uu_k+n(\psi_\al\nu^\al)_ku)
\end{aligned}
\end{equation}
we conclude that the expression is smooth in $x$ with uniformly bounded
$C^m$\nobreakdash- norms.

\cvm
The estimate \re{7.6} follows by induction.
\ep

\bl\lal{7.6}
The following estimates are valid
\begin{equation}\lae{7.10}
\norm{D\dot u}\le c \msp e^{-\ga t},
\end{equation}
\begin{equation}\lae{7.11}
\norm{\tfrac d{dt}F^{-1}}\le c\msp  F^{-1},
\end{equation}
and
\begin{equation}\lae{7.12}
\abs{\dot{\tilde v}}+\abs{\Ddot{\tilde v}}+\norm{D\dot{\tilde v}}\le c\msp
e^{-2\ga t}.
\end{equation}
Moreover,  $\dot{\tilde v}\msp e^{2\ga t}$ and $\ddot{\tilde v}\msp e^{2\ga t}$
converge, if
$t$ goes to infinity.
\el

\bp
\cq{\re{7.10}}\q The estimate follows immediately from
\begin{equation}
\dot u=\frac{\tilde v}F,
\end{equation}
in view of \rc{7.4}.

\cvm
\cq{\re{7.11}}\q Differentiating with respect to $t$ we obtain
\begin{equation}
\tfrac d{dt}F^{-1}=-F^{-2}(\dot H-n\dot{\tilde v}f'-n\tilde vf''\dot u+n\tfrac
d{dt}(\psi_\al\nu^\al))
\end{equation}
and the result follows from \re{7.12} and the known estimates for $\abs{\dot u}$
and $F$.

\cvm
\cq{\re{7.12}}\q We differentiate the relation $\tilde v=\h_\al\nu^\al$ to get
\begin{equation}\lae{7.15}
\begin{aligned}
\dot{\tilde v}&=\h_{\al\bet}\nu^\al\dot x^\bet+\h_\al\dot\nu^\al\\[\cma]
&=-\h_{\al\bet}\nu^\al\nu^\bet F^{-1}+(F^{-1})_ku^k
\end{aligned}
\end{equation}
yielding the estimate for $\abs{\dot{\tilde v}}$, in view of \rc{7.4} and the decay of
$\h_{\al\bet}$.

\cvm
Differentiating \re{7.15} covariantly with respect to $x$ we infer the estimate for
$\norm{D\dot{\tilde v}}$, while the estimate for $\norm{D\dot{\tilde v}}$ can be
 deduced after differentiating \re{7.15} covariantly with respect to $t$, in view
of \re{7.10}.

\cvm
The convergence of $\dot{\tilde v}\msp e^{2\ga t}$ and $\ddot{\tilde v}\msp
e^{2\ga t}$ can be easily verified.
\ep

Finally, let us estimate $\Ddot h^j_i$ and $\Ddot{\tilde h}^j_i$.

\bl\lal{7.7}
$\Ddot h^j_i$ and $\Ddot{\tilde h}^j_i$ decay like $e^{-\ga t}$.
\el

\bp
The estimate for $\Ddot h^j_i$ follows immediately by differentiating equation
\re{5.17} covariantly with respect to $t$ and by applying the above lemmata as well
as \rt{6.3}.

\cvm
Observing the remarks at the beginning of \rs{6} about rearranging crucial terms in
\re{5.17}, \cf equations \re{6.4} and \re{6.5}, we further conclude 
\begin{equation}
\norm{\Ddot {\tilde h}^j_i}\le c\msp e^{-\ga t}.\qedhere
\end{equation}
\ep

Using the same argument as in the proof of \rc{7.3} we infer

\bc\lac{7.8}
The tensor $\ddot h^j_i\msp e^{\ga t}$ converges, if $t$ tends to infinity.
\ec

The claims in \rt{0.3} are now almost all proved with the exception of two. In order to
prove the remaining claims we need

\bl
The function $\f=e^{\tilde\ga f}u^{-1}$ converges to $-\tilde\ga\sqrt m$ in
$C^\un(\so)$, if
$t$ tends to infinity.
\el

\bp
$\f$ converges to $-\tilde\ga\sqrt m$ in view of \re{3.5}. Hence, we only have to
show that
\begin{equation}
\norm{D^m\f}\le c_m\qq\A\,m\in\N^*,
\end{equation}
which will be achieved by induction.

\cvm
We have
\begin{equation}
\begin{aligned}
\f_i&=\tilde\ga e^{\tilde\ga f}f'u_iu^{-1}-e^{\tilde\ga f} u^{-2}u_i\\[\cma]
&=\f(\tilde\ga f'u-1)u^{-1}u_i.
\end{aligned}
\end{equation}

Now, we observe that
\begin{equation}
u^{-1}u=\tilde u^{-1}\tilde u_i
\end{equation}
and $f'u$ have uniformly bounded $C^m$\nobreak- norms in view of \rt{3.6},
\rl{7.1} and \rl{7.4}.

The proof of the lemma is then completed by a simple induction argument.
\eps

\bl
Let $(\breve g_{ij})$ be the induced metric of the leaves of the inverse mean
curvature flow, then the rescaled metric
\begin{equation}
e^{\tfrac2n t}\breve g_{ij}
\end{equation}
converges in $C^\un(\so)$ to
\begin{equation}
(\tilde\ga m)^{\frac1{\tilde\ga}} (-\tilde u)^\frac2{\tilde\ga}\bar\s_{ij},
\end{equation}
where we are slightly ambiguous by using the same symbol to denote $\tilde
u(t,\cdot)$ and
$\lim\tilde u(t,\cdot)$.
\el

\bp
There holds
\begin{equation}
\breve g_{ij}=e^{2f} e^{2\psi}(-u_iu_j+\s_{ij}(u,x)).
\end{equation}

Thus, it suffices to prove that
\begin{equation}\lae{7.24}
e^{2f}e^{\frac2n t}\ra (\tilde\ga m)^\frac1{\tilde\ga} (-\tilde u)^\frac2{\tilde\ga}
\end{equation}
in $C^\un(\so)$. But this evident in view of the preceding lemma, since
\begin{equation}
e^{2f}e^{\frac2n t}=(-e^{\tilde\ga f}u^{-1})^\frac2{\tilde\ga}(-\tilde
u)^\frac2{\tilde\ga}.\qedhere
\end{equation}
\ep

Finally, let us prove that the leaves $M(t)$ of the IMCF get more umbilical, if $t$
tends to infinity. Denote by $\breve h_{ij}, \breve\nu$, etc., the geometric
quantities of the hypersurfaces 
$M(t)$ with respect to the original metric $(\breve g_{\al\bet})$ in $N$,  then
\begin{equation}
e^{\tilde\psi}\breve h^j_i=h^j_i+\tilde\psi_\al\nu^\al\de^j_i,
\end{equation}
and hence,
\begin{equation}
\breve H^{-1}\abs{\breve h^j_i-\tfrac1n \breve H \de^j_i}=F^{-1}\abs{h^j_i-\tfrac1n H
\de^j_i}
\le c e^{-2\ga t}.
\end{equation}

\cvm
In  case $n+\om-4>0$, we even get a better estimate, namely,
\begin{equation}
\begin{aligned}
\abs{\breve h^j_i-\tfrac1n \breve H \de^j_i}&=e^{-\psi} e^{- f}e^{-\frac1n
t}\abs{h^j_i-\tfrac1n H\de^j_i}e^{\ga t} e^{(\frac1n -\ga)t}\\[\cma]
&\le c e^{-\frac1{2n}(n+\om-4)t},
\end{aligned}
\end{equation}
in view of \re{7.24}.

\section{Transition from big crunch to big bang}\las{8}

We shall define a new spacetime $\hat N$ by reflection and time reversal such that
the IMCF in the old spacetime transforms to an IMCF in the new one.

\cvm
By switching the light cone we obtain a new spacetime $\hat N$. The flow equation
in $N$ is independent of the time orientation, and we can write it as
\begin{equation}
\dot x=-\breve H^{-1}\breve\nu=-(-\breve H)^{-1}(-\breve\nu)\equiv -\hat
H^{-1}\hat \nu,
\end{equation}
where the normal vector $\hat \nu=-\breve\nu$ is past directed in $\hat N$ and the
mean curvature $\hat H=-\breve H$ negative.

Introducing a new time function $\hat x^0=-x^0$ and formally new coordinates
$(\hat x^\al)$ by setting
\begin{equation}
\hat x^0=-x^0,\q\hat x^i=x^i,
\end{equation}
we define a spacetime $\hat N$ having the same metric as $N$---only expressed in
the new coordinate system---such that the flow equation has the form
\begin{equation}\lae{8.3}
\dot{\hat x}=-\hat H^{-1}\hat \nu,
\end{equation}
where $M(t)=\graph \hat u(t)$, $\hat u=-u$, and 
\begin{equation}
(\hat\nu^\al)=-\tilde ve^{-\tilde \psi}(1,\hat u^i)
\end{equation}
in the new coordinates, since
\begin{equation}
\hat\nu^0=-\breve \nu^0\pde{\hat x^0}{x^0}=\nu^0
\end{equation}
and
\begin{equation}
\hat\nu^i=-\breve\nu^i.
\end{equation}

\cvm
The singularity in $\hat x^0=0$ is now a past singularity, and can be referred to as a
big bang singularity.

\cvm
The union $N\uu\hat N$ is a smooth manifold, topologically a product
$(-a,a)\ti\so$---we are well aware that formally the singularity $\{0\}\ti\so$ is not
part of the union; equipped with the respective metrics and time orientation it is a
spacetime which has a (metric) singularity in
$x^0=0$. The time function
\begin{equation}\lae{8.7}
\hat x^0=
\begin{cases}
\hp{-}x^0, &\tup{in } N,\\
-x^0, &\tup{in } \hat N,
\end{cases}
\end{equation}
is smooth across the singularity and future directed.

\cvm
$N\uu\hat N$ can be regarded as a \tit{cyclic universe} with a contracting part
$N=\{\hat x^0<0\}$ and an expanding part $\hat N=\{\hat x^0>0\}$ which are
joined at the singularity $\{\hat x^0=0\}$, \cf \cite{khoury:cyclic, steinhardt:cyclic} 
for similar ideas.

\cvm
We shall show that the inverse mean curvature flow, properly rescaled, defines a
natural $C^3$\nobreak- diffeomorphism across the singularity and with respect to
this diffeomorphism we speak of a transition from big crunch to big bang.

\cvm
Using the time function in \re{8.7} the inverse mean curvature flows in $N$ and
$\hat N$ can be uniformly expressed in the form
\begin{equation}\lae{8.8}
\dot{\hat x}=-\hat H^{-1}\hat\nu,
\end{equation}
where \re{8.8} represents the original flow in $N$, if $\hat x^0<0$, and the flow in
\re{8.3}, if $\hat x^0>0$.

\cvm
Let us now introduce a new flow parameter
\begin{equation}
s=
\begin{cases}
-\ga^{-1}e^{-\ga t},& \tup{for the flow in } N,\\
\hp{-}\ga ^{-1}e^{-\ga t},& \tup{for the flow in } \hat N,
\end{cases}
\end{equation}
and define the flow $y=y(s)$ by $y(s)=\hat x(t)$. $y=y(s,\x)$ is then defined in
$[-\ga^{-1},\ga^{-1}]\times \so$, smooth in $\{s\ne 0\}$, and satisfies the
evolution equation
\begin{equation}\lae{8.10}
y'\equiv \tfrac d{ds}y=
\begin{cases}
-\hat H^{-1}\hat\nu \msp e^{\ga t}, & s<0,\\
\hp{-}\hat H^{-1}\hat\nu \msp e^{\ga t},& s>0.
\end{cases}
\end{equation}

\bt
The flow $y=y(s,\x)$  is of class $C^3$ in $(-\ga^{-1},\ga^{-1})\times
\so$ and defines a natural diffeomorphism across the singularity. The flow parameter
$s$ can be used as a new time function.
\et

The flow $y$ is certainly continuous across the singularity, and also future directed,
i.e., it runs into the singularity, if $s<0$, and moves away from it, if $s>0$.

The continuous differentiability of $y=y(s,\x)$ with respect $s$ and $\x$ up to  order
three will be proved in a series of lemmata.

As in the previous sections we again view the hypersurfaces as embeddings with
respect to the ambient metric
\begin{equation}
d\bar s^2=-(dx^0)^2+\s_{ij}(x^0,x)dx^idx^j.
\end{equation}
The flow equation for $s<0$ can therefore be written as
\begin{equation}\lae{8.12}
y'=-F^{-1}\nu e^{\ga t}.
\end{equation}

\bl
$y$ is of class $C^1$ in $(-\ga^{-1},\ga^{-1})\ti \so$.
\el

\bp
Here, as in the proofs to come, we have to show that $y'$ and $y_i$ are continuous
in $\{0\}\ti\so$.

\cvm
Now, we have
\begin{equation}
y^0(s)=x^0(t),\q \hp{-} y^i(s)=x^i(t)\qq\A\,s<0,
\end{equation}
and
\begin{equation}
y^0(s)=-x^0(t),\q y^i(s)=x^i(t)\qq\A\,s>0,
\end{equation}
hence $y'$ is continuous across the singularity if and only if
\begin{equation}\lae{8.15}
\lim_{s\ua 0}\tfrac d{ds}y^0=\lim_{s\da 0}\tfrac d{ds}y^0,
\end{equation}
and
\begin{equation}\lae{8.16}
\lim_{s\ua 0}\tfrac d{ds}y^i=-\lim_{s\da 0}\tfrac d{ds}y^i.
\end{equation}

\cvm
Furthermore, we have to show that
\begin{equation}
\lim_{s\ua 0}y^0_i=0
\end{equation}
and
\begin{equation}
\lim_{s\ua 0}y^j_i=\lim_{s\da 0}y^j_i.
\end{equation}

The last two relations are obviously valid.

To verify \re{8.15} and \re{8.16} we observe

\br
The limit relations for $\spd {D^my}{\pde{}{x^0}}$ and $\spd{D^my}{\pde{}{x^i}}$,
where $D^my$ stands for covariant derivatives of order $m$ of $y$ with respect to
$s$ or $\x^i$, are identical to those for $\spd{D^my}\nu$ and $\spd{D^my}{x_i}$,
because $\nu$ converges to $-\pde{}{x^0}$, if $s\ua 0$.
\er

Thus, in view of \re{8.10} and \re{8.12}, it suffices to prove the convergence of
$Fe^{-\ga t}$, if $t$ goes to infinity. But this has already been shown in the proof of
\rc{5.4}, \cf equation \re{5.30}.
\ep

Let us examine the second derivatives.

\bl
$y$ is of class $C^2$ in $(-\ga^{-1},\ga^{-1})\ti\so$.
\el

\bp
\cq{$y'_i$}\q The normal component of $y_i'$ has to converge and the tangential
components have to converge to zero.

\cvm
We may only consider the behaviour for $s<0$. Then
\begin{equation}
y'=-F^{-1}e^{\ga t}\nu
\end{equation}
and
\begin{equation}
y'_i=F^{-2}F_ie^{\ga t}\nu -F^{-1}e^{\ga t} \nu_i
\end{equation}

The normal component is therefore equal to
\begin{equation}
F^{-2}e^{\ga t}(H_i-n\tilde v_if'-n\tilde vf''u_i+n\psi_{\al\bet}\nu^\al x^\bet_i,
+n\psi_\al \nu^\al_i)
\end{equation}
which converges to
\begin{equation}
\lim -F^{-2}e^{2\ga t} n f'' u^2\tilde u_i \tilde u^{-1}=\ga \tilde u \tilde u_i.
\end{equation}

The tangential components are equal to
\begin{equation}
-F^{-1}e^{\ga t} h^k_i,
\end{equation}
which converge to zero.

\cvm
\cq{$y_{ij}$}\q The Gau{\ss} formula yields
\begin{equation}
y_{ij}=h_{ij}\nu,
\end{equation}
which converges to zero as it should.

\cvm
\cq{$y''$}\q Here, the normal component has to converge to zero, while the
tangential ones have to converge.

\cvm
We get for $s<0$
\begin{equation}\lae{8.25}
\begin{aligned}
y''&=-\tfrac D{dt}(F^{-1}\nu)e^{2\ga t}-F^{-1}\nu \ga e^{2\ga t}\\[\cma]
&= -F^{-1}\dot\nu e^{2\ga t}+F^{-2}\nu \dot F e^{2\ga t} -F^{-1}\nu\ga
e^{2\ga t}.
\end{aligned}
\end{equation}

The normal component is equal to
\begin{equation}
F^{-2}e^{2\ga t}(\dot H-n\dot{\tilde v}f'-n\tilde vf''\dot u +n\psi_{\al\bet}\nu^\al
\dot x^\bet +n\psi_\al\dot\nu^\al -\ga F).
\end{equation}

$F^{-2}e^{2\ga t}$ converges, all other terms converge to zero with the possible
exception of
\begin{equation}
-n\tilde vf''\dot u-\ga F=-F^{-1}n(\tilde v^2 f'' +\tfrac1n \ga F^2)
\end{equation}
which however converges to zero too, in view of \re{4.37} and the estimate for
$\abs{H}$.

\cvm
The tangential components are equal to
\begin{equation}
\begin{aligned}
F^{-1}D_i(F^{-1})e^{2\ga t}&=-F^{-3}e^{2\ga t}(H_i-n\tilde v_i f' -n\tilde v
f''u_i\\[\cma]
&\qq\qq+n\psi_{\al\bet}\nu^\al x^\bet_i+n\psi_\al \nu^\al_i),
\end{aligned}
\end{equation}
which converge to
\begin{equation}
\lim F^{-3}e^{3\ga t} n\tilde v f''u^2\tilde u_i\tilde u^{-2}.\qedhere
\end{equation}
\ep

\bl
$y$ is of class $C^3$ in $(-\ga^{-1},\ga^{-1})\ti\so$.
\el

\bp
\cq{$y_{ijk}$}\q Now, the normal component has to converge to zero, while the
tangential ones have to converge. Again we look at $s<0$ and get
\begin{equation}\lae{8.30}
y_{ij}=h_{ij}\nu,
\end{equation}
\begin{equation}
y_{ijk}=h_{ijk}\nu+h_{ij}\nu_k.
\end{equation}
Hence, $y_{ijk}$ converges to zero.

\cvm
\cq{$y'_{ij}$}\q The normal component has to converge, while the tangential ones
should converge to zero.

\cvm
Using the Ricci identities it can be easily checked that, instead of $y'_{ij}$, we may
look at $\tfrac D{ds}(y_{ij})$, since
\begin{equation}\lae{8.32}
\riema 0i0j\ra 0,
\end{equation}
\cf \rl{4.4}.

From \re{8.30} we deduce
\begin{equation}
\tfrac D{ds}y_{ij}=\dot h_{ij}\nu e^{\ga t} +h_{ij}\dot\nu e^{\ga t},
\end{equation}
and conclude further that the normal component converges, in view of \rc{7.3}, and
the tangential ones converge to zero, since $\dot\nu$ vanishes in the limit.

\cvm
\cq{$y''_i$}\q The normal component has to converge to zero and the tangential
ones have to converge.

\cvm
From \re{8.25} we infer
\begin{equation}\lae{8.34}
\begin{aligned}
y''&=-F^{-3}e^{2\ga t}(H^k-n\tilde v^kf'-nf''u^k+n(\psi_\al\nu^\al)^k)x_k\\[\cma]
&\hp{=}+F^{-2}e^{2\ga t}(\dot H-n\dot{\tilde v}f'+n\tfrac
D{dt}(\psi_\al\nu^\al))\nu\\[\cma]
&\hp{=}+F^{-3}e^{2\ga t}(-n\tilde v^2[f''+\tilde\ga \abs{f'}^2]-\ga
[H^2+n^2(\psi_\al\nu^\al)^2-2nHf'\tilde v\\[\cma]
&\hp{=}\qq +2nH\psi_\al\nu^\al -2n^2 f'\tilde v \psi_\al\nu^\al])\nu
\end{aligned}
\end{equation}
and thus
\begin{equation}
\begin{aligned}
y''_i&=-(F^{-3}e^{2\ga t}(H^k-n\tilde
v^kf'-nf''u^k+(n\psi_\al\nu^\al)^k))_ix_k\\[\cma]
&\hp{=}-F^{-3}e^{2\ga t}(H^k-n\tilde
v^kf'-nf''u^k+(n\psi_\al\nu^\al)^k)h_{ik}\nu\\[\cma]
&\hp{=}+(F^{-2}e^{2\ga t}(\dot H-n\dot{\tilde v}f'+\tfrac
D{dt}(n\psi_\al\nu^\al)))_i\nu\\[\cma]
&\hp{=}+F^{-2}e^{2\ga t}(\dot H-n\dot{\tilde v}f'+\tfrac
D{dt}(n\psi_\al\nu^\al))\nu_i\\[\cma]
&\hp{=}+(F^{-3}e^{2\ga t}(-n\tilde v^2[f''+\tilde\ga \abs{f'}^2]-\ga
[H^2+(n\psi_\al\nu^\al)^2-2nHf'\tilde v\\[\cma]
&\hp{=}\qq +2Hn\psi_\al\nu^\al -2n f'\tilde v n\psi_\al\nu^\al]))_i\nu\\[\cma]
&\hp{=}+F^{-3}e^{2\ga t}(-n\tilde v^2[f''+\tilde\ga \abs{f'}^2]-\ga
[H^2+(n\psi_\al\nu^\al)^2-2nHf'\tilde v\\[\cma]
&\hp{=}\qq +2Hn\psi_\al\nu^\al -2n f'\tilde v n\psi_\al\nu^\al])\nu_i.
\end{aligned}
\end{equation}

\cvm
Therefore, the normal component converges to zero, while the tangential ones
converge.

\cvm
\cq{$y'''$}\q The normal component has to converge, while the tangential ones have
to converge to zero.

\cvm
Differentiating the equation \re{8.34} we get
\begin{equation}
\begin{aligned}
&y'''=3F^{-4}e^{3\ga t}\dot F(H^k-n\tilde v^kf'-nf''u^k+(n\psi_\al\nu^\al)^k)x_k
\\
&\hp{=}-2\ga F^{-3}e^{3\ga t}(H^k-n\tilde
v^kf'-nf''u^k+(n\psi_\al\nu^\al)^k)x_k\\
&\hp{=}-F^{-3}e^{3\ga t}\tfrac D{dt}(H^k-n\tilde
v^kf'-nf''u^k+(n\psi_\al\nu^\al)^k)x_k\\
&\hp{=}-F^{-3}e^{3\ga t}(H^k-n\tilde
v^kf'-nf''u^k+(n\psi_\al\nu^\al)^k)\dot x_k\\
 &\hp{=}-2F^{-3}e^{3\ga t}\dot
F(\dot H-n\dot{\tilde v}f'+\tfrac D{dt}(n\psi_\al\nu^\al))\nu\\
&\hp{=}+2\ga F^{-2}e^{3\ga t}
(\dot H-n\dot{\tilde v}f'+\tfrac D{dt}(n\psi_\al\nu^\al))\nu\\
&\hp{=}+ F^{-2}e^{3\ga t}
\tfrac D{dt}(\dot H-n\dot{\tilde v}f'+\tfrac D{dt}(n\psi_\al\nu^\al))\nu\\
&\hp{=}+ F^{-2}e^{3\ga t}
(\dot H-n\dot{\tilde v}f'+\tfrac D{dt}(n\psi_\al\nu^\al))\dot\nu\\
&\hp{=}-3F^{-4}e^{3\ga t}\dot F (-n\tilde v^2[f''+\tilde\ga \abs{f'}^2]-\ga
[H^2+(n\psi_\al\nu^\al)^2-2nHf'\tilde v\\[\cma]
&\hp{=}\qq +2Hn\psi_\al\nu^\al -2n f'\tilde v n\psi_\al\nu^\al])\nu\\[\cma]
&\hp{=}+2\ga F^{-3}e^{3\ga t}(-n\tilde v^2[f''+\tilde\ga \abs{f'}^2]-\ga
[H^2+(n\psi_\al\nu^\al)^2-2nHf'\tilde v\\[\cma]
&\hp{=}\qq +2Hn\psi_\al\nu^\al -2n f'\tilde v n\psi_\al\nu^\al])\nu\\[\cma]
&\hp{=}+F^{-3}e^{3\ga t}\tfrac D{dt}(-n\tilde v^2[f''+\tilde\ga \abs{f'}^2]-\ga
[H^2+(n\psi_\al\nu^\al)^2-2nHf'\tilde v\\[\cma]
&\hp{=}\qq +2Hn\psi_\al\nu^\al -2n f'\tilde v n\psi_\al\nu^\al])\nu\\[\cma]
&\hp{=}+F^{-3}e^{3\ga t}(-n\tilde v^2[f''+\tilde\ga \abs{f'}^2]-\ga
[H^2+(n\psi_\al\nu^\al)^2-2nHf'\tilde v\\[\cma]
&\hp{=}\qq +2Hn\psi_\al\nu^\al -2n f'\tilde v n\psi_\al\nu^\al])\dot\nu\\[\cma]
\end{aligned}
\end{equation}

\cvm
Observing that
\begin{equation}
\dot x_k=F^{-2}F_k\nu-F^{-1}\nu_k
\end{equation}
and
\begin{equation}
\dot u_k=F^{-1}\tilde v_k-F^{-2}\tilde vF_k
\end{equation}
and taking the results of \rl{7.6}, \rl{7.7}, and \rc{7.8} into account we conclude that
the normal component converges.

\cvm
The tangential components contain the following crucial terms
\begin{equation}
\begin{aligned}
&3F^{-4}e^{3\ga t}n^2\tilde v^2\abs{f''}^2u^k\dot u+2\ga F^{-3}e^{3\ga t}n\tilde
vf''u^k\\[\cma]
&+F^{-3}e^{3\ga t}n\tilde vf'''u^k\dot u+F^{-5}e^{3\ga t}n^2\tilde
v^3\abs{f''}^2u^k,
\end{aligned}
\end{equation}
which can be rearranged to yield
\begin{equation}
F^{-5}e^{3\ga t}n\tilde v u^k(4f''(f''+\tilde \ga
\abs{f'}^2)-f'(f''+\tilde\ga\abs{f'}^2)').
\end{equation}
Hence, the tangential components tend to zero.

\cvm
The remaining mixed derivatives of $y$, which are obtained by commuting the order
of differentiation in the derivatives we already treated, are also continuous across the
singularity in view of the Ricci identities and \re{8.32}.
\ep

\section{ARW spaces and the Einstein equations}

Let $N$ be a cosmological spacetime such that the metric has the form as specified
in \rd{0.1}, though, with regard to $f$, we only assume at the moment that $f$ is
smooth and satisfies
\begin{equation}\lae{9.1}
\lim_{\tau\ra b}f(\tau)=-\un
\end{equation}
and
\begin{equation}\lae{9.2}
\lim_{\tau\ra b}f'= -\un.
\end{equation}

The conformal metric
\begin{equation}\lae{9.3}
d\bar s^2=e^{2\psi}(-(dx^0)^2+\s_{ij}(x^0,x)dx^idx^j)
\end{equation}
should satisfy the conditions in \rd{0.1}, and, in addition, the partial derivatives of
$\psi$ as well as the second fundamental form  of the coordinate slices
$\{x^0=\const\}$ and its derivatives should be integrable over the range
$[a,b)$ of
$x^0$.

In contrast to the previous sections we suppose that the Einstein equations are valid
\begin{equation}
G_{\al\bet}=\ka T_{\al\bet},
\end{equation}
where $\ka$ is a positive constant, and the stress-energy tensor is asymptotically
equal to that of a perfect fluid.

\bd
Let $x^0$ be a time function such the preceding assumptions are satisfied. A
symmetric, divergent free tensor $(T_{\al\bet})$ is said to be asymptotically equal to
that of a perfect fluid with respect to the future, if the mixed tensor $(T^\al_\bet)$
splits in the form
\begin{equation}
T^\al_\bet=\bar T^\al_\bet +\hat T^\al_\bet,
\end{equation}
where $(\bar T^\al_\bet)$ is the stress-energy tensor of a perfect fluid, i.e.,
\begin{equation}
\bar T^0_0=-\rho,\q T^\al_i=\de^\al_ip;
\end{equation}
$0\le \rho$ is the \tit{density} and $p$ is the \tit{pressure}, and $(\hat
T^\al_\bet)$ as well as its partial derivatives of arbitrary order are supposed to
vanish, if $x^0$ tends to $b$, and they should be integrable over the range $[a,b)$ of
$x^0$. Moreover, $\hat T^\al_\bet f'$ should vanish and  be integrable as well.
\ed

\cvm
Let us assume an \tit{equation of state}
\begin{equation}
p=\tfrac\om{n}\rho
\end{equation}
holds, where $\om\in\R[]$ is a constant such that
\begin{equation}
n+\om-2>0.
\end{equation}

We shall now show that, because of the Einstein equations, $f$ has to satisfy the
conditions stated in \rd{0.1}, even slightly stronger ones.

First, we prove
\bl\lal{9.2}
There exist $\tau_0$ and $c>0$ such that
\begin{equation}
\rho(\tau,x)\ge c>0\qq\A\,\tau\ge \tau_0,\;  \A x\,\in\so.
\end{equation}
Moreover,
\begin{equation}
\lim_{\tau\ra b}\rho=\un.
\end{equation}
\el

\bp
We use the Einstein equations
\begin{equation}
G_{00}=\ka T_{00}
\end{equation}
to conclude
\begin{equation}
\tfrac12n(n-1)\abs{f'}^2+\tfrac12\bar R+\e=\ka \rho e^{2\tilde\psi}+\ka \hat
T_{00},
\end{equation}
where we recall that $\bar R$ is the scalar curvature of the metric in \re{9.3},
and where
$\e$ represents terms that converge to zero, if $\tau$ tends to $b$, or equivalently,
\begin{equation}\lae{9.13}
\tfrac12n(n-1)\abs{f'}^2e^{-2\tilde\psi}+\tfrac12\bar R e^{-2\tilde\psi}+\e\msp
e^{-2\tilde\psi}=\ka
\rho +\ka\hat T^0_0.
\end{equation}
Hence, we have
\begin{equation}
\ka\rho\sim\tfrac12 n(n-1)\abs{f'}^2e^{-2f},
\end{equation}
which proves the result, in view of \re{9.1} and \re{9.2}.
\ep

\bl
Let $\tilde\ga=\tfrac12(n+\om-2)$, then there exists a constant $m>0$ such that
\begin{equation}\lae{9.15}
\lim_{\tau\ra b}\abs{f'}^2e^{2\tilde\ga f}=m
\end{equation}
and
\begin{equation}\lae{9.16}
\abs{D^mf}\le c_m \abs{f'}^m\qq\A\,m\in\N.
\end{equation}
Furthermore, the limit metric $(\bar\s_{ij})$ must have constant scalar curvature
$R$.  The function 
\begin{equation}
\f=f''+\tilde\ga\abs{f'}^2,
\end{equation}
converges to
\begin{equation}\lae{9.18}
\lim_{\tau\ra b}\f=-\tfrac\ga{n-1}R,
\end{equation}
where $\ga=\tfrac1n{\tilde\ga}$, and in addition
\begin{equation}\lae{9.19}
\lim_{\tau\ra b}D^m\f=0\qq\A\,m\in\N^*.
\end{equation}
\el

\bp
\cq{\re{9.15}}\q Since $(T_{\al\bet})$ is divergent free, we deduce
\begin{equation}
\begin{aligned}
0=T^\ga_{0;\ga}&=T^\ga_{0,\ga}+\breve\C^\ga_{\ga\al}T^\al_0-\breve
\C^\al_{0\ga}T^\ga_\al\\[\cma]
&=-\dot\rho-(n+\om)\dot{\tilde\psi}\rho
-\tfrac12\dot\s^i_i(1+\tfrac\om{n})\rho+\mc C,
\end{aligned}
\end{equation}
where $\mc C$ tends to zero and is integrable over the range  $[a,b)$ of $x^0$.

\cvm
In view of \rl{9.2} we deduce
\begin{equation}
\tfrac d{dt}\log\rho=-(n+\om)\dot{\tilde\psi}+\mc C,
\end{equation}
where we still use the same symbol $\mc C$, and hence, for fixed $x$,
\begin{equation}
\rho(\tau,x)e^{(n+\om)\tilde\psi(\tau,x)}=\rho(\tau',x)e^{(n+\om)\tilde\psi(\tau',x)}
e^{\int_{\tau'}^\tau \mc C}.
\end{equation}

Thus, we conclude, first, that $\rho(\tau,x)e^{(n+\om)\tilde\psi(\tau,x)}$ is
uniformly bounded, and then, that it converges to a positive function, if $\tau$ tends
to $b$.

At the moment the limit can depend on the spatial variables $x$, but we shall see
immediately that it is a constant.

\cvm
Now, multiplying equation \re{9.13} with $e^{(n+\om)\tilde\psi}$ we deduce
\begin{equation}
\lim\abs{f'}^2e^{(n+\om-2)f}=\tfrac2{n(n-1)}\ka\lim \rho e^{(n+\om)f},
\end{equation}
i.e., the limit on the left-hand side exists, and the limit on the right-hand side is a
constant.

\cvm
\cq{\re{9.16}}\q We consider the contracted version of the Einstein equations
\begin{equation}
G^\al_\al=\ka
T^\al_\al
\end{equation}
and infer with the help of equation \re{9.13}
\begin{equation}
\begin{aligned}
\breve R&=\tfrac2{n-1}\ka\rho(1-\om)+\mc
C\\[\cma]
&=n\abs{f'}^2e^{-2\tilde\psi}(1-\om)+\tfrac1{n-1}\bar Re^{-2\tilde\psi}(1-\om)+\e
e^{-2\tilde\psi}+\mc C,
\end{aligned}
\end{equation}
and we further conclude
\begin{equation}\lae{9.26}
\tfrac\ga{n-1}\bar R+f''+\tfrac12(n+\om-2)\abs{f'}^2=\e+\mc C e^{2\tilde\psi}.
\end{equation}

The estimate in \re{9.16} now follows immediately by induction.

\cvm
\cq{\re{9.18} and \re{9.19}}\q One easily checks that
\begin{equation}
\lim_{\tau\ra b}\bar R=R,
\end{equation}
where $R$ is the scalar curvature of $(\bar\s_{ij})$.

The relation \re{9.26} implies that $\f$ is uniformly Lipschitz continuous and
bounded, hence there exists a sequence $\tau_k\ra b$ such that $\f(\tau_k)$
converges, from which we deduce that $R$ has to be constant. Therefore,
$\f=f''+\tilde\ga\abs{f'}^2$ converges.

\cvm
Moreover, after having established the relation \re{9.15}, we can apply the result of
\rl{3.1}, i.e., $b$ is finite, and without loss of generality we may assume that $b=0$,
which in turn allows us to conclude that derivatives of arbitrary order of the
right-hand side of \re{9.26} tend to zero in the limit, \cf \rl{4.4}.

This completes the proof of the lemma.
\ep

\bibliographystyle{amsplain}
\providecommand{\bysame}{\leavevmode\hbox to3em{\hrulefill}\thinspace}
\providecommand{\MR}{\relax\ifhmode\unskip\space\fi MR }
\providecommand{\MRhref}[2]{%
  \href{http://www.ams.org/mathscinet-getitem?mr=#1}{#2}
}
\providecommand{\href}[2]{#2}



\end{document}